\documentclass[12pt]{amsart}

\usepackage{amsmath,amsthm,amscd,amsfonts,amssymb,epic,eepic,bbm,stmaryrd}
\usepackage[pagebackref,colorlinks=true,linkcolor=blue,citecolor=blue]{hyperref}

\allowdisplaybreaks
\newtheorem{thm}{Theorem}[section]
\newtheorem{cor}[thm]{Corollary}

\newtheorem{lem}[thm]{Lemma}

\newtheorem{prop}[thm]{Proposition}
\theoremstyle{remark}
\newtheorem*{rem}{Remark}

\newcounter{remarkscounter}

\numberwithin{equation}{section}
\newcommand{\A}{\mathbb{A}}
\newcommand{\GL}{\mathrm{GL}}
\newcommand{\SL}{\mathrm{SL}}
\newcommand{\ZZ}{\mathbb{Z}}

\newcommand{\Gal}{\mathrm{Gal}}

\newcommand{\QQ}{\mathbb{Q}}

\newcommand{\lto}{\longrightarrow}

\newcommand{\OO}{\mathcal{O}}
\newcommand{\CC}{\mathbb{C}}
\newcommand{\RR}{\mathbb{R}}
\newcommand{\GG}{\mathbb{G}}

\newcommand{\gl}{\mathfrak{gl}}
\newcommand{\quash}[1]{}

\theoremstyle{definition}

\renewcommand{\bar}{\overline}

\numberwithin{equation}{subsection}

\newcommand{\one}{\mathbbm{1}}

\begin{document}
\title[The Rankin-Selberg monoid and a nonabelian trace formula]{A summation formula for the Rankin-Selberg monoid and a nonabelian trace formula}
\author{Jayce R. Getz}
\address{Department of Mathematics\\
Duke University\\
Durham, NC 27708}
\email{jgetz@math.duke.edu}


\thanks{The author is thankful for partial support provided by NSF grant DMS-1405708.  Any opinions, findings, and conclusions or recommendations expressed in this material are those of the author and do not necessarily reflect the views of the National Science Foundation.}

\maketitle

\begin{abstract}

Let $F$ be a number field and let $\mathbb{A}_F$ be its ring of adeles.  Let $B$ be a quaternion algebra over $F$ and let $\nu:B \lto F$ be the reduced norm.
Consider the reductive monoid $M$ over $F$ whose points in an $F$-algebra $R$ are given by
\begin{align*}
M(R):=\{(\gamma_1,\gamma_2) \in (B \otimes_F R)^{2}:\nu (\gamma_1)=\nu(\gamma_2)\}.
\end{align*}
Motivated by an influential conjecture of Braverman and Kazhdan we prove a summation formula analogous to the Poisson summation formula for certain spaces of functions on the monoid.  As an application, we define new zeta integrals for the Rankin-Selberg $L$-function and prove their basic properties.  We also use the formula to prove a nonabelian twisted trace formula, that is, a trace formula whose spectral side is given in terms of automorphic representations of the unit group of $M$  that are isomorphic (up to a twist by a character) to their conjugates under a simple nonabelian Galois group.
\end{abstract}

\tableofcontents

\section{Introduction}

Let $F$ be a number field, let $\A_F$ be its ring of adeles and let $f \in \mathcal{S}(\gl_n(\A_F))$, the Schwartz space of $\gl_n(\A_F)$. 
The Poisson summation formula on $\gl_n(F)$ is
\begin{align} \label{PS}
\sum_{\gamma \in \gl_n(F)} f(\gamma)
=\sum_{\gamma \in \gl_n(F)}\widehat{f}\left(\gamma\right)
\end{align}
where $\widehat{f}$ is the Fourier transform of $f$ with respect to an additive character.   Godement and Jacquet used this formula to prove the functional equation for the standard $L$-functions of automorphic representations of $\GL_n(\A_F)$ \cite{GodementJacquetBook}.  Motivated by Godement and Jacquet's work, Braverman and Kahzdan \cite{BK-lifting} have conjectured that for every representation 
$$
\rho:{}^L G^{\circ} \to \GL(V)
$$ 
of the neutral component of the $L$-group ${}^LG$ of a connected reductive group $G$ there should be a corresponding $\rho$-Fourier transform and $\rho$-Poisson summation formula that in turn imply the analytic continuation and functional equation of the Langlands $L$-function $L(s,\pi,\rho)$. 
 L.~Lafforgue 
has a related program \cite{LafforgueJJM}. We also note that Ng\^o has advocated investigating Braverman and Kazhdan's conjecture using the trace formula.    For more information on the Fourier transforms see \cite{BouthierNgoSakellaridis,Cheng:Ngo:BK,WWLiSat,WWLi:Zeta} in the nonarchimedean case and \cite{GetzBK,WWLi:Towards} for the archimedean case.  See also \cite{BK:normalized,Getz:Liu:BK,Pollack:GJ,Shahidi:FT} for more information about cases where much more is known about Braverman and Kazhdan's conjectures.

In this paper we investigate this circle of ideas when the monoid $\gl_n$ is replaced by a monoid related to the Rankin-Selberg convolution on $\GL_2 \times \GL_2$.  We then apply the work to produce a nonabelian trace formula isolating automorphic representations that are invariant (up to an abelian twist) to their Galois conjugates under a simple nonabelian group.

Let $B$ be a simple algebra of degree $4$ over $F$; thus $B$ is either a quaternion algebra over $F$ or $M_2(F)$.   Let $\OO_B <B$ be the order of \eqref{OB}.
Let
$$
\nu: B \lto F
$$ 
be the reduced norm.  
Consider the group scheme whose points in an $\OO_F$-algebra $R$ are given by 
$$
G(R)=\{(g_1,g_2) \in ((\OO_B \otimes_{\OO_F}R)^\times)^2: \nu (g_1)=\nu (g_2)\}.
$$
The neutral component of the $L$-group ${}^LG:={}^LG_{F}$  is
$$
{}^{L}G^{\circ}=\GL_2(\CC) \times \GL_2(\CC)/\{(zI,z^{-1}I):z \in \CC^\times\}
$$
 and hence the tensor product induces a representation
$$
\rho:{}^LG^{\circ} \lto \GL_{4}(\CC).
$$
Now $G$ is the group of units of the monoid scheme whose points in an $\OO_F$-algebra $R$ are given by
\begin{align} \label{M:def}
M(R):=\{(X_1,X_2) \in (\OO_B \otimes_{\OO_F}R)^{2}:\nu( X_1)=\nu( X_2) \}.
\end{align}
We refer to this as the \textbf{Rankin-Selberg monoid}. 
In the next subsection we define local spaces of functions related to the monoid $M$ and define a pair of transforms using them.
We use these functions to state an asymptotic formula for sums over $M(F)$ in \S \ref{ssec:asymptotic}.  In \S \ref{ssec:zeta} we then define zeta integrals for  $\rho$ analogous to those introduced by Godement and Jacquet for the standard representation.  We will prove their analytic continuation using the asymptotic formula mentioned above.  The nonabelian trace formula will be discussed in \S \ref{ssec:ntf}.

Throughout this paper we use the following notation:
For an $F$-algebra $R$ and $T=(T_1,T_2) \in B_R^{2}$ we let
\begin{align*}
p(T):&=\nu (T_1)-\nu( T_2),\\
\mathrm{tr}\, T:&=\mathrm{tr}(T_1+T_2), \nonumber\\
\nu (T):&=\nu (T_1)\nu( T_2),\\ \nonumber
\omega(T):&=\nu(T_1).
\end{align*}  
Here $\mathrm{tr}$ denotes the reduced trace and $\nu$ denotes the reduced norm.  Also, if $\pi$ is a representation of $G(\A_F)$ or $G(F_v)$ and $s \in \CC$ then
$$
\pi_s(g):=\pi(g)|\nu(g)|^{s/2}=\pi(g)|\omega(g)|^s
$$
where $|\cdot|$ is either the norm on $\A_F^\times$ or $F_v^\times$ depending on the context.

\subsection{Local spaces of functions}

\label{ssec:funcs}

In this subsection we use local notation, fixing a place $v$ of $F$ and writing, e.g.~$F:=F_v$.
Let $\psi:F \to \CC^\times$ be a nontrivial character.  If $W$ is an $F$-vector space we let $\mathcal{S}(W)$ be the usual Schwartz space (it is just $C_c^\infty(W)$ in the nonarchimedean case).  Throughout this work $\Phi$ will denote a function in $\mathcal{S}(W)$; the space $W$ will depend on the context.

For $\Phi \in \mathcal{S}(F \times F \times B_F^2)$ define
\begin{align} \label{IPhi}
\mathcal{I}(\Phi,s):B_F^{2} \lto \CC
\end{align}
via 
\begin{align} \label{IPhichi1}
\mathcal{I}(\Phi,s)(\gamma):&=\int_{F^{\times}}\left(
\int_{B^{2}_{F}}\Phi\left(\frac{p(T)}{t},t,T\right)\psi\left(\frac{\mathrm{tr}\,\gamma T}{t}\right) dT\right)|t|^sdt^\times.
\end{align}
Here and below $dT$ (resp.~$dt$) denotes the additive Haar measure on $B_F$ (resp.~$F$), both normalized so that the Fourier transform is self-dual with respect to $\psi$, and 
$$
dt^\times=\zeta(1)\frac{dt}{|t|}
$$ 
with $\zeta(s)$ the usual local Euler factor of the Dedekind zeta function
(see \S \ref{ssec-notation}).    The integrals depend on the choice of $\psi$, but we will not encode this into the notation.   
It is clear that these integrals converge absolutely for $\mathrm{Re}(s)>0$, and they in fact admit a meromorphic continuation to $\mathrm{Re}(s)>-4$ to a function that is holomorphic at $s=-3$ (see Theorem \ref{thm:loc:cont}).
   We set
\begin{align} \label{IPhi}
\mathcal{I}(\Phi)(\gamma):=\zeta(1)^{-1}\mathcal{I}(\Phi,-3)(\gamma).
\end{align}

If $F$ is nonarchimedean we let $\varpi$ be a choice of uniformizer and $q:=|\OO_F/\varpi|$. 
  For  $g_1,g_2 \in G(F)$ let
\begin{align} \label{basic:func}
\one_{\rho}(g_1,g_2)=\sum_{k=0}^\infty q^{2k}\one_{\OO_B}(\varpi^{-k}g_1,\varpi^{-k}g_2)
\end{align}
where $\OO_B < B$ is a maximal order.  If $B$ is split then this is the \textbf{basic function} attached to the representation $\rho$. 
By this we mean that for an unramified representation $\pi$ of $G(F)$ and $\mathrm{Re}(s)$ sufficiently large one has that 
$
\pi_{s+1}(\one_{\rho})$ projects the space of $\pi$ onto the line fixed by $G(\OO_F)$ and acts via the scalar $L(s,\pi,\rho)$ on that line (see Lemma \ref{lem:sat:comp}).

It turns out that if $q$ is sufficiently large, $\psi$ is unramified, and the maximal order $\OO_B$ is chosen 
appropriately then
\begin{align} \label{Ibasic}
\mathcal{I}(\one_{\OO_F^2 \times \OO_B^{2}})|_{G(F)}=\frac{\one_{\rho}}{\zeta(2)}
\end{align} for a suitable choice of Haar measures (see Theorem \ref{thm:unram}). 
Let 
\begin{align} \label{I:space}
\mathcal{I}(G(F)):=\{\mathcal{I}(\Phi):\Phi \in \mathcal{S}(F \times F \times B_F^{ 2})\}
\end{align}
where $\mathcal{I}(\Phi)$ is defined as in \eqref{IPhi}.
\begin{rem}One might ask why this space of functions is not called $\mathcal{I}(B_F^{2})$, since its elements are defined on all of $B_F^{2}$.  It is largely a matter of emphasis; we are most interested in the restriction of these functions to $M(F)$.  We do not denote the space by $\mathcal{I}(M(F))$ because the value of the functions at $(0,0)$ is only given via analytic continuation, and hence if we view these as functions on $M(F)$ it is difficult to make sense of the equality \eqref{IG} below.
\end{rem}

We are not sure whether $\mathcal{I}(G(F))$ is the Schwartz space for which Braverman and Kazhdan are searching, but it is certainly closely related to it.  
 Regardless of whether this is the ``actual'' Schwarz space attached to $\rho$ the global results of this paper (i.e.~theorems \ref{thm:zeta} and \ref{thm:ntf}) show that it is of interest.
For $(x,y,T) \in F \times F \times B_F^{ 2}$ let 
\begin{align}
\Phi^{\mathrm{sw}}(x,y,T)=\Phi(y,x,T)
\end{align}
(the $\mathrm{sw}$ is for ``switch'').
We then have a pair of transforms
\begin{align}
\mathcal{I}(G(F)) \longleftarrow &\mathcal{S}(F \times F \times B_F^{2})\lto \mathcal{I}(G(F))\\
\mathcal{I}(\Phi)  \longmapsfrom &\quad \quad \quad  \Phi \quad \quad \quad \quad \longmapsto \mathcal{I}(\Phi^{\mathrm{sw}}). \nonumber
\end{align}
The functions $\mathcal{I}(\Phi)$ and $\mathcal{I}(\Phi^{\mathrm{sw}})$ behave somewhat like a Fourier transform pair (see Theorem \ref{thm:zeta}).

\subsection{The asymptotic formula} \label{ssec:asymptotic}
We now revert to global notation.  
Let $S$ be a finite set of places of $F$ including the
 infinite places.
We define adelic space of functions $\mathcal{I}(G(\A_F))$ via 
\begin{align} \label{IG}
\mathcal{I}(G(\A_F))=\otimes_v' \mathcal{I}(G(F_v))
\end{align}
where the restricted direct product is taken with respect to the basic functions $\one_{\rho,v}$.  
For $\Phi \in \mathcal{S}(\A_F^2 \times B_{\A_F}^{ 2})$ and $s \in \CC$ we denote by $\mathcal{I}(\Phi,s)$ the obvious adelic analogue of the local integral  \eqref{IPhichi1}  of the previous section.  It converges absolutely for $\mathrm{Re}(s)>1$, and admits a meromorphic continuation to the half plane $\mathrm{Re}(s)>-4$ (see Theorem \ref{thm:glob:cont}).  
For $g \in G(\A_F)$ we set
\begin{align} \label{IPhi:adelic}
\mathcal{I}(\Phi)(g):=\frac{\mathrm{Res}_{s=-3}\mathcal{I}(\Phi,s)(g)}{\mathrm{Res}_{s=1}\Lambda(s)}  \quad \textrm{  and } \quad
\mathcal{I}_{-1}(\Phi)(0):=\frac{\mathrm{Res}_{s=-1}\mathcal{I}(\Phi,s)(0)}{\mathrm{Res}_{s=1}\Lambda(s)}.
\end{align} 
Here $\Lambda(s)$ is the completed Dedekind zeta function of $F$.  

 The function $\mathcal{I}(\Phi)$ is in fact defined on  $B^{ 2}$ but vanishes if $\gamma \not \in M(F)$ (see Theorem \ref{thm:glob:cont}).  Moreover $\mathcal{I}(\Phi)$, as a function of $G(\A_F)$, is in $\mathcal{I}(G(\A_F))$.  Indeed, if $\Phi=\Phi_S \one_{(\widehat{\OO}_F^S)^2 \times (\widehat{\OO}_B^{S})^2}$ for a sufficiently large finite set of places $S$ including the infinite places 
 then
$$
\mathcal{I}(\Phi)=\frac{\mathcal{I}(\Phi_S)}{\zeta^S(2)}\prod_{v \not \in S}\one_{\rho,v}
$$ 
by Theorem \ref{thm:unram}.  Here and in global settings below $\zeta$ denotes the Dedekind $\zeta$-function of $F$.

For $\Phi_0 \in \mathcal{S}(\A_F^2)$ let 
$$
\mathcal{F}_2(\Phi_0)(x,y)=\int_{\A_F} \Phi_0(x,t)\psi(yt)dt
$$ 
be the Fourier transform of $\Phi_0$ in the second variable.
For our later use, we say that a Schwartz function $\Phi \in \mathcal{S}(\A_F^2 \times B_{\A_F}^2)$ \textbf{satisfies the standard assumptions} if 
there exists
$$
(\Phi_0,f) \in  \mathcal{S}(\A_F^2) \times \mathcal{S}(B_{\A_F}^{2})
$$
such that $\Phi=\Phi_0 \otimes f$ and 
\begin{enumerate}
\item $\Phi_0(t,0)=0$ for all $t \in \A_F$ 
\item $\mathcal{F}_2(\Phi_0)(0,0)=\Gamma_{F_\infty}(1)$
\end{enumerate}
where $\Gamma_{F_\infty}(s)$ is the factor of $\Lambda(s)$ at infinity.  Notice that if $\Phi$ satisfies the standard assumptions, then it is impossible for $\Phi^{\mathrm{sw}}$ to satisfy the standard assumptions because of (1) and (2).

If $X \in \RR_{>0}$ we denote by
\begin{align} \label{Delta:def}
\Delta(X) \in \A_F^{\times}
\end{align}
the idele that is $X^{[F:\QQ]^{-1}}$ at all places $v|\infty$ and $1$ elsewhere.  We also denote by $\Delta$ the isomorphism
\begin{align*}
\RR_{>0} &\tilde{\lto} A_G\\
X &\longmapsto \Delta(X)(I,I)
\end{align*}
where $A_G$ is the usual central subgroup of $G(F_\infty)$ (see \S \ref{ssec-notation}).  We endow $A_G$ with the Haar measure corresponding to $\frac{dx}{x}$ via this isomorphism.

The starting point of this paper is the following theorem, which will play the role of the Poisson summation formula in our setting:
\begin{thm} \label{thm:main:asymp} 
For any $\varepsilon>0$ the sum 
\begin{align} \label{basic:sum}
\sum_{\gamma \in M(F)}f\left( \frac{\gamma}{\Delta(\sqrt{X})} \right)
\end{align}
is equal to $O_{\varepsilon,f,\Phi}(X^{\tfrac{3}{2}+\varepsilon})$ plus
\begin{align*}
X^3\left(\mathcal{I}_{-1}(\Phi)(0)-\mathcal{I}_{-1}(\Phi^{\mathrm{sw}}) (0)\right)
+
X^{2}\sum_{\gamma \in M(F)}  \left(\mathcal{I}(\Phi)(\gamma)-\mathcal{I}(\Phi^{\mathrm{sw}})(\gamma)\right).
\end{align*}
Moreover
\begin{align*}
\sum_{\gamma \in M(F)}  \left|\mathcal{I}(\Phi)(\gamma)\right|+\left|\mathcal{I}(\Phi^{\mathrm{sw}})(\gamma)\right|<\infty.
\end{align*} 
\end{thm}

\noindent We will explain how to derive Theorem \ref{thm:main:asymp} from \cite[Theorem 1.1]{GetzQuad} in \S \ref{sec:translate} below.
 It is worth noting that the proof of Theorem \ref{thm:main:asymp} makes no use of automorphic representations.

\subsection{Zeta functions for Rankin-Selberg convolutions} \label{ssec:zeta}

 Using Theorem \ref{thm:main:asymp} we can execute the argument of Godement and Jacquet to obtain new zeta integrals for Rankin-Selberg convolutions as we now explain.

Let $\pi$ be an infinite-dimensional automorphic representation of $A_{G} \backslash G(\A_F)$ and let $\varphi$ be a smooth function in the space of $\pi$.    For $\Phi \in \mathcal{S}(\A_F^2 \times B_{\A_F}^2)$ let
\begin{align*}
Z(\mathcal{I}(\Phi),s,\varphi):=\int_{G(\A_F)} \mathcal{I}(\Phi)(g)|\omega(g)|^{s+1}\varphi(g)dg.
\end{align*}
Moreover, for $g \in G(\A_F)$ let
\begin{align}
g^{\vee}:=\omega(g)^{-1}g\quad \textrm{ and } \quad \varphi^{\vee}(g):=\varphi(g^{\vee}).
\end{align}
 We note that if $\pi$ is an irreducible automorphic representation of $G(\A_F)$  then the representation
$$
g \mapsto \pi(g^{\vee})
$$  
is isomorphic to the contragredient because $B$ is a simple algebra of rank $2$ (this is false in higher rank). 

Let $\Phi=\Phi_0 \otimes f$ satisfy the standard assumptions, and assume moreover that $f=f_S \one_{(\widehat{\OO}_B^S)^2}$ where $f_S \in C_c^\infty(G(F_S))$ for some finite set of places $S$ including the infinite places.  Throughout this paper we make the following assumption on $S$:
\begin{enumerate}
\item[\textbf{A($S$)}] The finite set $S$ contains all infinite places, all dyadic places, all places dividing the absolute different of $F$, the places dividing the elements $a,b$ defining the quaternion algebra $B$ and all places where $\psi$ is ramified. 
\end{enumerate}

We also assume that for $a \in A_{G}$ and $g \in G(\A_F)^1$  one has
\begin{align} \label{h:prod}
f_S(\Delta(a)g)=V(a)f_S^1(g)
\end{align}
for some $f_S^1 \in C_c^\infty(G(F_S)^1)$ and some $V \in C_c^\infty(\RR_{>0})$. Here the superscript $1$ groups are defined as in \S \ref{ssec-notation}.
Assuming that $f_S$ decomposes as above 
is convenient in the proof and is no loss of generality spectrally.  

The following is the first main theorem of this this paper:

\begin{thm} \label{thm:zeta} Assume $B$ is a division algebra.  The zeta functions $Z(\mathcal{I}(\Phi),s,\varphi)$ and $Z(\mathcal{I}(\Phi^{\mathrm{sw}}),s,\varphi^{\vee})$ admit meromorphic continuations to the plane and satisfy a functional equation
$$
Z(\mathcal{I}(\Phi),s,\varphi)=Z(\mathcal{I}(\Phi^{\mathrm{sw}}),1-s,\varphi^\vee).
$$
The function $Z(\mathcal{I}(\Phi),s,\varphi)$ is holomorphic except for a possible simple pole at $s=1$.
\end{thm}

Let $S$ be a sufficiently large finite set of places including the infinite and dyadic places.  By  \eqref{Ibasic} if $\varphi$ is spherical outside of $S$ then 
\begin{align*}
Z(\mathcal{I}(\Phi),s,\varphi)=Z(\mathcal{I}(\Phi_S),s,\varphi)\frac{L^S(s,\pi,\rho)}{\zeta^S(2)}
\end{align*}
where 
$$
Z(\mathcal{I}(\Phi_S),s,\varphi)=\int_{G(F_S)}\mathcal{I}(\Phi_S)(g_S)|\omega(g_S)|_S^{s+1}\varphi(g_SI^S)dg_S
$$
(here $I^S$ is the element of $G(\A_F^S)$ that is the identity outside of $S$).  Thus Theorem \ref{thm:zeta} implies a coarse version of the functional equation of the Rankin-Selberg 
$L$-function.  Our proof, moreover, is new. It uses the summation formula Theorem \ref{thm:main:asymp} in place of the Langlands-Shahidi method or Rankin-Selberg theory.   

We remark that in analogy with Tate zeta functions one might expect that there are two poles of the zeta integrals $Z(\mathcal{I}(\Phi),s,\varphi)$, one at $0$ and one at $1$.  
We expect that absence of the pole at $s=0$ can be traced either to our assumption that $f_S \in C_c^\infty(G(F_S))$ or assumption (1) in our standard assumptions on $\Phi$.  Making this precise would at very least  require reworking much of the paper \cite{GetzQuad}.

\subsection{A nonabelian trace formula}
\label{ssec:ntf}
In this subsection we explain the second main theorem of the paper.  It is a nonabelian twisted trace formula, that is, a twisted trace formula whose spectral side is given in terms of automorphic representations of $G(\A_F)$ that are isomorphic to their conjugates under a group of automorphisms of $F$ generated by a pair of elements.  Since any simple nonabelian group can be generated by two elements \cite[Corollary 8.3]{Guralnick_Malle_2012}, this is a quite general setup.  

Assume that there is a subfield $k \leq F$ such that $F/k$ is Galois with Galois group 
$$
\Gal(F/k)=\langle\iota, \tau \rangle;
$$ 
that is, $\Gal(F/k)$ is generated by two elements.  For example $\Gal(F/k)$ could be any simple nonabelian group \cite[Corollary 8.3]{Guralnick_Malle_2012}.  Assume moreover that $B_1$ is a division algebra over $k$ such that $B:=B_1 \otimes_kF$ is nonsplit (i.e.~again a division algebra).

Let $G_0$ be the connected reductive $k$-subgroup of $\mathrm{Res}_{F/k}G$ 
whose points in a $k$-algebra $R$ are given by
$$
G_0(R):=\{(g_1,g_2) \in G(R \otimes_kF)=(B^{2} \otimes_k R)^\times: \nu (g_1)=\nu( g_2) \in R^\times\}.
$$
This $k$-subgroup comes equipped with an action of $\Gal(F/k)^{2}$, and in particular we have an automorphism $\theta$ of $G_0$ given on points by 
\begin{align*}
\theta(g_1,g_2):=(\iota(g_1),\tau(g_2)).
\end{align*}
We have an action of $G_0$ on $\mathrm{Res}_{F/k}G$ via $\theta$ conjugation, given on points in an $F$-algebra $R$ by 
\begin{align*}
G_0(R) \times \mathrm{Res}_{F/k}G(R) &\lto \mathrm{Res}_{k/F}G(R)\\
(g_1,\gamma) &\longmapsto g_1 \gamma \theta(g_1)^{-1}.
\end{align*}
For $\gamma \in G(F)=\mathrm{Res}_{F/k}G(k)$ we let $G_{0\gamma}$ be the stabilizer of $\gamma$ under this action.  Since $B$ is a division algebra, a standard argument implies that $G_{0\gamma}$ is reductive and anisotropic modulo center.  For suitable test functions $f$ on $G(\A_F)$ we can then form twisted orbital integrals $\mathrm{TO}_{\gamma}(f)$ in the usual manner (see \eqref{TO}).

For an automorphic representation $\pi$ and $f \in C_c^\infty(G(\A_F))$ we introduce the following \textbf{nonabelian trace}:
\begin{align} \label{ntr}
\mathrm{ntr}\,\pi(f)=\frac{\mathrm{Res}_{s=1}L^S(s,\pi,\rho)}{\zeta^S(2)}\int_{[G_0]}K_{\pi|\omega|^2 (f)}(g,\theta(g))dg.
\end{align}
Here $[G_0]$ is the typical adelic quotient (see \eqref{quotient}).
It is clear that the integral here is absolutely convergent since $[G_0]$ is compact.  It is also not hard to see that it vanishes unless the $L$-packet of $\pi$ is stable under $\Gal(F/k)$ up to a twists by abelian characters; for a precise statement we refer to Lemma \ref{lem:ntr}. 

The following is the last main theorem of this paper.  It will be proven in \S \ref{sec:nabel}:

\begin{thm} \label{thm:ntf} 
One has
\begin{align*}
\sum_{\pi} \mathrm{ntr}\, \pi(f)=
&\mathrm{meas}([G_0])\left(\mathcal{I}(\Phi)(0)-\mathcal{I}(\Phi^{\mathrm{sw}})(0)\right)\\
&+\sum_{\gamma} \mathrm{meas}([G_{0\gamma}])\left(\mathrm{TO}_{\gamma}(\mathcal{I}(\Phi))-\mathrm{TO}_{\gamma}(\mathcal{I}(\Phi^{\mathrm{sw}}))\right)
\end{align*}
where the sum on $\pi$ is over isomorphism classes of automorphic representations of $A_G \backslash G(\A_F)$ and the sum on $\gamma$ is over a set of representatives for the orbits of $G_0(F)$ acting on $G(F)$ via $\theta$-conjugation. 
All of the sums here are absolutely convergent.
\end{thm}

Ultimately, one would like to compare this formula with an analogous one over $k$ and prove nonsolvable base change and descent for automorphic representations of inner forms of $\GL_2$.  More details on what one should expect spectrally from such a comparison are contained in \cite{GetzApproach}.  Unfortunately we have no paradigm available to inform such a comparison.  

 A smaller step towards this goal would be to use  Theorem \ref{thm:ntf} to study limit multiplicities of forms that are isomorphic to their $\Gal(F/k)$ conjugates as the analytic conductor of the corresponding automorphic representations goes to infinity in some fashion.  One expects that this is an easier problem, and it is still of great interest.

\subsection{Outline of the paper}

In \S \ref{sec:translate} we state the results from \cite{GetzQuad} that we require for this paper and explain how to translate them into the current setting.  We describe a little of the structure of $L$-packets for $G$ in \S \ref{sec:Lpackets}.  This is required for \S \ref{sec:spectral}, in which we give a spectral expansion of the sum \eqref{basic:sum}.  Here we make no use of Rankin-Selberg theory.  We then prove the analytic continuation and functional equations of our zeta functions $Z(\mathcal{I}(\Phi),s,\varphi)$ in \S \ref{sec:GJ}.  Again, we make no use of Rankin-Selberg theory, so this work provides a new proof of the meromorphic continuation of $L^S(s,\pi,\rho)$.

In \S \ref{sec:four:var} we allow ourselves to use the entirety of Rankin-Selberg theory and use it to prove an identity for what we call a four-variable kernel function.  It is an identity between a spectral side and a geometric side.  The key point here is that the spectral side is given as a sum over kernel functions attached to automorphic representations $\pi$ of $B_{\A_F}^\times$, but there are four variables of integration attached to each $\pi$, not just two as in the Selberg expansion of the standard automorphic kernel.  These extra variables are the key to proving Theorem \ref{thm:ntf}.  We carry out the argument in 
\S \ref{sec:nabel}.  

\section*{Acknowledgements}

The author thanks V.~Blomer, T.~Kaletha, L.~Pierce, A.
~Pollack, D.~Schindler and W.~Zhang for useful conversations. 
H.~Hahn read the entire paper several times as it was being written; the author truly appreciates her constant encouragement and help with editing.  He also thanks P.~E.~Herman for pointing out the reference \cite{DukeIwaniecConvolution}.

\section{Notation}
\label{ssec-notation}
Let $\psi:F \backslash \A_F \to \CC^\times$ denote a nontrivial additive character.  For $\Phi \in \mathcal{S}(B_{\A_F})$ we let
$$
\widehat{\Phi}(Y):=\int_{B_{\A_F}}\Phi(X)\psi(\mathrm{tr}( YX)) dX
$$
denote the Fourier transform of $\Phi$.  
We always normalize the Haar measure on $B_{\A_F}$ so that it is self-dual with respect to this transform.  
 The Poisson summation formula then takes the form
\begin{align*}
\sum_{\gamma \in B}\Phi(\gamma)=\sum_{\gamma \in B}\widehat{\Phi}(\gamma).
\end{align*}

It is convenient to use the following additional (fairly standard) notation.  For
an affine algebraic group $H$ over a number field $F$ let $A_H$ be the connected component in the real topology of the greatest $\QQ$-split torus in $\mathrm{Res}_{F/\QQ}Z_H$ ($Z_H$ being the center of $H$).  We then set 
\begin{align} \label{quotient}
[H]:=A_H H(F) \backslash H(\A_F).
\end{align}
We will only use this notation when $H$ is reductive, in which case this quotient has finite measure with respect to the left-invariant measure induced by a Haar measure on $H(\A_F)$.  One has the usual Harish-Chandra map 
\begin{align} \label{HC:map} \begin{split}
\mathrm{HC}_H:H(\A_F) &\lto \mathrm{Lie}\,A_H=\mathrm{Hom}_{\ZZ}(X(\mathrm{Res}_{F/\QQ}H),\RR)\\
x &\longmapsto (\chi \mapsto |\log (\chi(x))|) \end{split}
\end{align}
where $X(\mathrm{Res}_{F/\QQ}H)$ is the character group of $\mathrm{Res}_{F/\QQ}H$.
We let
\begin{align} \label{HC:group}
H(\A_F)^1:=\mathrm{ker}\, \mathrm{HC}_H
\end{align}
and 
let $H(F_S)^1$ be the kernel of the composite of the natural inclusion $H(F_S) \hookrightarrow H(\A_F)$ and $\mathrm{HC}_H$.  

Finally if $S$ is a finite set of places of $F$ including the infinite places then 
$\OO_F^S \subset F$ is the ring of $S-\infty$-integers, and 
$$
\widehat{\OO}^S_F:=\prod_{v \not \in S} \OO_{F_v},
$$
which is isomorphic to the profinite completion of $\OO_F^S$.  In particular $\widehat{\OO}^S_F \cap F=\OO^S_F$.

\section{Preliminaries} \label{sec:translate}

\subsection{The proof of Theorem \ref{thm:main:asymp}}

The setting of \cite{GetzQuad} involves a vector space $V=\GG_a^n$ and a nondegenerate quadratic form $Q$ on $V$.  Given these data and 
$\Phi \in \mathcal{S}(V(\A_F))$, $\chi \in \widehat{[\GG_m]}$ one defines 
\begin{align} \label{IPhi:old}
\mathcal{I}(\Phi,\chi_s)(\xi):=\int_{\A_F^\times \times V(\A_F)}\Phi\left(\frac{Q(w)}{t},t,w \right)\psi\left( \frac{\langle \xi,w\rangle}{t}\right)dw\chi_s(t)dt^\times
\end{align}
Here $\langle, \rangle:\GG_a^n \times \GG_a^n \to \GG_a$ is the usual product:
\begin{align}
\langle a,b \rangle=\sum_{i=1}^n a_ib_i.
\end{align}
For each place $v$ of $F$ we have the analogous local integrals.

We explain how the setting of the current paper fits into this setup.  Here and below we will freely use basic facts about quaternion algebras.  A new and enjoyable reference is \cite{Voight:QA}.
Let $\left(\frac{a,b}{F} \right)$ be the quaternion algebra over $F$ 
consisting of 
$$
\{x_1+x_2i+x_3j+x_4k: x_i \in F\}
$$
where $i,j,k$ are subject to the usual relations: $ij=k=-ji$ and $i^2=a$, $j^2=b$.
Fix $a,b$ such that $B \cong \left(\frac{a,b}{F} \right)$; we use this isomorphism to identify $B$ and $\left(\frac{a,b}{F} \right)$.  We assume without loss of generality that $a,b \in \OO_F$.  Let
\begin{align} \label{OB}
\OO_B:=\{x_1+x_2i+x_3j+x_4k: x_i \in \OO_F\}.
\end{align}
Throughout this paper if $v$ is a finite place of $F$ we let
$$
\OO_{B_{F_v}}:=\OO_B \otimes_{\OO_F} \OO_{F_v}
$$
and $\widehat{\OO}_B^S:=\prod_{v \not \in S} \OO_{B_{F_v}}$.  

We use a bar to denote the canonical involution
\begin{align*}
\left( \frac{a,b}{F} \right) &\lto \left(\frac{a,b}{F} \right)\\
x_1+x_2i+x_3j+x_4k &\longmapsto x_1-x_2i-x_3j-x_4k.
\end{align*}

We identify $B$ and $F^4$ via the (ordered) basis $1,i,j,k$.
Writing $x=x_1+x_2i+x_3j+x_4k$, $y=y_1+y_2i+y_3j+y_4k$ we have
\begin{align*}
\nu(x)&=x_1^2-ax_2^2-bx_3^2+abx_4^2,\\
\tfrac{1}{2}\mathrm{tr}(xy)&=x_1y_1+ax_2y_2+bx_3y_3-abx_4y_4.
\end{align*}

In this paper our vector space $V$ is $B \oplus B$, identified with $(F^4)^2=F^8$ and our quadratic form is $p(T)=\nu(T_1)-\nu(T_2)$.
With respect to the basis given by 
$$
((1,0),(i,0),(j,0),(k,0),(0,1),(0,i),(0,j),(0,k))
$$ 
the matrix $J$ of this form is diagonal with entries 
$$
(1,-a,-b,ab,-1,a,b,-ab).
$$
Let $\alpha:B \to B$ be the $F$-vector space isomorphism given by
$$
\alpha(x_1+x_2i+x_3j+x_4k)=2x_1+2ax_2i+2bx_3j-2abx_4k.
$$
Extend $\alpha$ to an isomorphism $\alpha:B^{2} \to B^{2}$ by letting it act on each factor separately.

Then the relationship between \eqref{IPhi} and \eqref{IPhi:old} in the case at hand is
\begin{align} \label{int:rel}
\mathcal{I}(\Phi,s)(\xi)=\mathcal{I}(\Phi,1_s)(\alpha(\xi)).
\end{align}
where $1_s(x):=|x|^s$.  The local analogue of this identity also holds.
Most of the results of \cite{GetzQuad} regarding the analytic properties of $\mathcal{I}(\Phi,s)$ involve
a character $\mathcal{G}$ (see \cite[Lemma 3.1]{GetzQuad}).    
Since  $(-1)^{8/2}\det J=a^4b^4$  is a square in $F$  this character is trivial in the case at hand \cite[Lemma 6.2]{GetzQuad}). 

In a moment we require the following lemma:
\begin{lem} \label{lem:prod:rule}
For $x \in F_\infty^n$ and $v|\infty$ let $|x|_v:=\max\{|x_i|_v:1 \leq i \leq n\}$.  Let $A>0$, $N>0$, $\beta \in \OO_F \cap F^\times$ be given.  If $\alpha \in \beta^{-1}\OO_F^n - 0$ then for $N' \in \ZZ_{>0}$ large enough in a sense depending on $A$, $N$ and $\beta$ one has
\begin{align*}
\prod_{v|\infty}\left(\max(|\alpha|_v,1)^{-N'}\min(|\alpha|_v,1)^{-A}\right)\ll_{A,N,\beta} \prod_{v|\infty}\max(|\alpha|_v,1)^{-N}.
\end{align*}

\end{lem}

\begin{proof}
We note first that it suffices to prove the lemma when $\beta=1$.  Indeed, if the lemma is true for $\beta=1$, then to deduce the general case one uses the fact that for $\beta \in F^\times$ one has
\begin{align*}
\max\left(\frac{|\alpha|_v}{|\beta|_v},1\right) &=\frac{1}{|\beta|_v}\max(|\alpha|_v,|\beta|_v) \asymp_{\beta} \max(|\alpha|_v,1)\\
\min\left(\frac{|\alpha|_v}{|\beta|_v},1\right) &=\frac{1}{|\beta|_v}\min(|\alpha|_v,|\beta|_v) \asymp_{\beta} \min(|\alpha|_v,1).
\end{align*}

Let $\alpha \in \OO_F^n- \{0\}$.  Then 
$$
\prod_{v|\infty}|\alpha|_v \geq \max_i |\alpha_i|_\infty \geq 1.
$$ 
 In 
particular 
if $|\alpha|_{v_1}<1$ for some $v_1|\infty$, there is another $v_2|\infty$ such that 
$$
|\alpha|_{v_2}\geq
\left(\frac{1}{|\alpha|_{v_1}}\right)^{(i-1)^{-1}}
$$
where $i$ is the number of infinite places of $F$. 
We can therefore take $N'=A(i-1)^2+N$.
\end{proof}

 For each place $v$ of $F$ we let 
\begin{align} \label{maxnorm}
|(x_1+x_2i+x_3j+x_4k,y_1+y_2i+y_3j+y_4k)|_v:=\max_{i,j}(|x_i|,|y_j|).
\end{align}

Using the lemma we can give good analytic control on sums involving $\mathcal{I}(\Phi)$:

\begin{thm} \label{thm:adelic:bound} Let $\Omega \subset G(\A_F) \times G(\A_F)$ be a compact set.  For $\gamma \in B^2$, $a \in A_G$ with $|\omega(a)| \geq 1$ and $(g_1,g_2) \in \Omega$ one has
\begin{align*}
|\mathcal{I}(\Phi)(ag_1^{-1}\gamma g_2)| \ll_{\Phi,N,\Omega} |\omega(a)|^{-N/2}\prod_{v|\infty}\max(|\gamma|,1)^{-N}.
\end{align*}
\end{thm}

\begin{proof}
We may assume $\gamma \neq 0$.  By \cite[Theorems 4.1, 5.2]{GetzQuad} there is an $A>0$ such that for any $N' \in \ZZ_{\geq 0}$ one has 
\begin{align} \label{withgs}
|\mathcal{I}(\Phi)(ag_1^{-1}\gamma g_2)| \ll_{\Phi,N'} \left(\prod_{v|\infty}(|ag_1^{-1}\gamma g_2|_v,1)^{-N'}\min(|ag_1^{-1}\gamma g_2|_v,1)^{-A}\right) \prod_{v \nmid \infty} \max(1,|g_1^{-1}\gamma g_2|_v^{-3}).
\end{align}
Moreover, $\mathcal{I}(\Phi)(ag_1^{-1}\gamma g_2)$ vanishes for $\gamma \not \in \beta^{-1}\OO_B$ for some $\beta \in F^\times$ (see \cite[Lemma 5.1]{GetzQuad}). 
 Here $\beta$ depends on $\Omega$, but not on $g_1,g_2$.  
Since $|ag_1^{-1}\gamma g_2|_v \asymp_{\Omega} |a\gamma |_v$ for all $v$ and 
$|ag_1^{-1}\gamma g_2|_v=|\gamma|_v$ for $v$ outside a finite set of places depending on $\Omega$ we have
\begin{align} \label{wogs} \begin{split}
|\mathcal{I}(\Phi)(ag_1^{-1}\gamma g_2)| &\ll_{\Phi,N',\Omega} \left(\prod_{v|\infty}\max(|a\gamma|_v,1)^{-N'}\min(|a\gamma|_v,1)^{-A}\right) \prod_{v \nmid \infty} \max(1,|\gamma|_v^{-3})\\
&\ll_{\beta,N'} |\omega(a)|^{-N'/2}\left(\prod_{v|\infty}\max(|\gamma|_v,1)^{-N'}\min(|\gamma|_v,1)^{-A}\right) \prod_{v \nmid \infty} \max(1,|\gamma|_v^{-3}).\end{split}
\end{align}
 Using Lemma \ref{lem:prod:rule} we obtain the bound
\begin{align*}
 |\mathcal{I}(\Phi)(ag_1^{-1}\gamma g_2)| \ll_{\Phi,\Omega,N}|\omega(a)|^{-N/2} \left(\prod_{v|\infty}\max(|\gamma|,1)^{-N}\right) \prod_{v \nmid \infty} \max(1,|\gamma|_v^{-3})
 \end{align*}
for any $N \in \ZZ_{\geq 0}$.  The factor involving the
 places $v \nmid \infty$ can then be absorbed at the expense of increasing $N$ and the implicit constant.
\end{proof}

The following theorem is a combination of the $\chi=1$ cases of
 \cite[Theorems 4.1, 4.2, 5.2]{GetzQuad} in the current setting: 

\begin{thm} \label{thm:loc:cont} For each place $v$ and all $\Phi \in \mathcal{S}(F_v^2 \oplus B_{F_v}^2)$ the integrals $\mathcal{I}(\Phi_v,s)(\gamma)$ admit meromorphic continuations to $\mathrm{Re}(s)>-4$ that are holomorphic at $s=-3$. \qed
\end{thm}

For $x$ in a vector space let
\begin{align*}
\delta_x:=\begin{cases} 1 & \textrm{ if } x=0 \\ 0 &\textrm{otherwise.}\end{cases}
\end{align*}
The following theorem is \cite[Theorem 3.2]{GetzQuad} in the current setting:
\begin{thm} \label{thm:glob:cont} For $\Phi \in \mathcal{S}(\A_F^2 \oplus B_{\A_F}^2)$ and $\gamma \in B^2$ the adelic integral $\mathcal{I}(\Phi,s)(\gamma)$ admits a meromorphic continuation to $\mathrm{Re}(s)>-4$.  Moreover, 
$$
\frac{\mathcal{I}(\Phi,s)}{s^{\delta_{\gamma}}(s+1)^{\delta_\gamma}(s+3)^{\delta_{p(\gamma)}}}
$$
is holomorphic for $\mathrm{Re}(s)>-4$. \qed
\end{thm}

The following is a special case of \cite[Theorem 6.4]{GetzQuad}:

\begin{thm} \label{thm:unram} Let $v$ be a finite place of $F$ that does not divide $2ab$.
 Assume in addition that $\psi_v$ is unramified and $F_v$ is absolutely unramified.  Then
\begin{align*}
\mathcal{I}(\one_{\OO_F^2 \oplus \OO_B^2})|_{G(F)}=\frac{\one_{\rho}}{\zeta(2)}.
\end{align*} \qed
\end{thm}

\begin{proof}[Proof of Theorem \ref{thm:main:asymp}]
Let $p^{\vee}$ be the quadratic form on $B^2$ whose matrix is $J^{-1}$.   Recall that we have just proven that the character $\mathcal{G}$ of \cite[Lemma 3.1]{GetzQuad} is trivial.  With this in mind we apply \cite[Theorem 3.3]{GetzQuad} in the setting above to see that 
\begin{align*}
\sum_{\substack{\gamma \in B^2\\ p(\gamma)=0}}f\left( \frac{\gamma}{\Delta(\sqrt{X})} \right)
\end{align*}
is equal to $O_{\varepsilon,f,\Phi}(X^{\tfrac{3}{2}+\varepsilon})$ plus
\begin{align*}
X^3\left(\mathcal{I}_{-1}(\Phi)(0)-\mathcal{I}_{-1}(\Phi^{\mathrm{sw}}) (0)\right)
+
X^{2}\sum_{\substack{\gamma \in B^2\\ p^{\vee}(\gamma)=0}}  \left(\mathcal{I}(\Phi)(\alpha^{-1} (\gamma))-\mathcal{I}(\Phi^{\mathrm{sw}})(\alpha^{-1} (\gamma))\right).
\end{align*}
Moreover using Theorem \ref{thm:adelic:bound} we see that
\begin{align*}
\sum_{\substack{\gamma \in B^2\\ p^{\vee}(\gamma)=0}}  \left|\mathcal{I}(\Phi)(\alpha^{-1}(\gamma))\right|+\left|\mathcal{I}(\Phi^{\mathrm{sw}})(\alpha^{-1}(\gamma))\right|<\infty.
\end{align*} 
We now note that $2p^\vee(\gamma)=p(\alpha^{-1}(\gamma))$.  Thus taking a change of variables $\gamma \mapsto \alpha(\gamma)$ we deduce the theorem.  
\end{proof}

\section{$L$-packets for $G$} \label{sec:Lpackets}

In this section we recall some results on $L$-packets for $G$ that will be necessary in the subsequent sections.  For $\OO_F$-algebras $R$ let
$$
\widetilde{G}(R):=\{g \in ((\OO_B \otimes_{\OO_F}R)^\times)^2\}.
$$
Thus (if we base change to $F$) $G$ is a subgroup of $\widetilde{G}$ and they share the same derived group, an inner form of $\SL_2^{2}$.  We let
\begin{align} \label{Res:L:map}
\mathrm{Res}: {}^L\widetilde{G}\lto{}^LG 
\end{align}
denote the natural quotient map.

\subsection{Local $L$-packets}

In this subsection we fix a place $v$ of $F$ and omit it from notation, writing $F:=F_v$.
Let $\widetilde{\pi}$ be an irreducible admissible representation of $\widetilde{G}(F)$.  The restriction $\widetilde{\pi}|_{G(F)}$ is a finite direct sum of irreducible admissible representations of $G(F)$.  It defines a set $\mathrm{Res}(\widetilde{\pi})$ of isomorphism classes of irreducible admissible representations of $G(F)$:
\begin{align*}
\mathrm{Res}(\widetilde{\pi}):=\{\pi: \pi \textrm{ is isomorphic to an irreducible subrepresentation of } \widetilde{\pi}|_{G(F)}\}.
\end{align*}
An $L$-packet of irreducible admissible representations of $G(F)$ is a set of representations of this form.  The set of $L$-packets of $G(F)$ partition the set of irreducible representations of $G(F)$.  For all of these facts see \cite[Chapter 2]{Hiraga:Saito} or \cite{Tadic:SLn} and the references therein.  We note that when $\widetilde{G} \cong \mathrm{GL}_2^2$ then each representation in $\mathrm{Res}(\widetilde{\pi})$ occurs with multiplicity $1$ in $\widetilde{\pi}|_{G(F)}$ \cite[Lemma 2.6]{Labesse:Langlands}.
If $\pi \in \mathrm{Res}(\widetilde{\pi})$ then we say that $\pi$ is a \textbf{transfer} of $\widetilde{\pi}$.  

\begin{lem} \label{lem:spherical}
Assume $F$ is nonarchimedean and let $\widetilde{\pi}$ be an unramified representation of $\widetilde{G}(F)$.  Then
there is a unique unramified irreducible representation in $\mathrm{Res}(\widetilde{\pi})$. 
\end{lem}
\begin{proof}
It is clear that there is an unramified irreducible subrepresentation in $\mathrm{Res}(\widetilde{\pi})$.
Note that $G^{\mathrm{der}} \cong \mathrm{SL}_2^2$.  
By \cite[Proposition 3.2.4]{Lansky:Raghuram} there is a unique unramified subrepresentation of the restriction $\widetilde{\pi}|_{G^{\mathrm{der}}(F)}$.  This is actually stronger than the first assertion of the lemma.
\end{proof}

This notion of transfer is compatible with the unramified local Langlands correspondence in the following sense:

\begin{lem} \label{lem:compat} Assume that $F$ is nonarchimedean.  Let $\widetilde{\pi}$ be an unramified representation of $\widetilde{G}(F)$.  The representation $\pi$ of $G(F)$ obtained from this representation using the $L$-map $\mathrm{Res}:{}^L\widetilde{G} \to {}^LG$ is the unique unramified representation in $\mathrm{Res}(\widetilde{\pi})$.
\end{lem}

\begin{proof}
Let $T_2 \leq \GL_{2\CC}$ be the maximal torus of diagonal matrices.  Then $ T_2 \times T_2 \leq {}^LG^\circ$ is a maximal torus.  Let $\widehat{Z} \leq T_2 \times T_2$ be the torus whose points in a $\CC$-algebra $R$ are 
$$
\widehat{Z}(R):=\{(zI,z^{-1}I):z \in R^\times\}.
$$ 
Let $\widehat{T}:=(T_2 \times T_2)/\widehat{Z}$.  
One has a diagram
\begin{align*}
\begin{CD}
C_c^\infty(G(F)//G(\OO_F)) @>{\mathcal{S}}>> \CC[\widehat{T}]^{W(\widehat{T},{}^LG^{\circ})}\\
@V{\mathrm{Res}^*}VV @V{\mathrm{Res}^*}VV\\
C_c^\infty(\widetilde{G}(F)//\widetilde{G}(\OO_F)) @>{\mathcal{S}}>> \CC[T_2 \times T_2]^{W(T_2 \times T_2,{}^L\widetilde{G}^{\circ})}
\end{CD}
\end{align*}
where $\mathcal{S}$ is the Satake isomorphism and the map on the right is that induced by $\mathrm{Res}:{}^L\widetilde{G}^\circ \to {}^L G^{\circ}$.  The map on the left is just the unique map that makes the diagram commutative.    If we view a function in $\CC[\widehat{T}]$ as a function on $T_2 \times T_2$ that is invariant under the diagonal action of $\widehat{Z}$ then this map is just given by the natural inclusion.  In particular, it is injective, and its image is precisely the elements of 
$\CC[T_2 \times T_2]^{W(T_2 \times T_2,{}^L\widetilde{G}^{\circ})}$ invariant under the action of $\widehat{Z}$.  Temporarily denote this algebra by $\mathcal{A}$; thus $\mathrm{Res}$ induces an isomorphism
$$
\mathrm{Res}^*:\CC[\widehat{T}]^{W(\widehat{T},{}^LG^{\circ})} \tilde{\lto} \mathcal{A}.
$$
 Let us describe $\mathcal{A}$.

Let $V_{\mathrm{st}}$ be the standard representation of $\GL_{2\CC}$.  Denote by $V_1$ and $V_2$ be the representations obtained by composing the natural projections ${}^L\widetilde{G}^\circ \lto \GL_{2\CC}$ with the homomorphism defining $V_{\mathrm{st}}$.  A basis of $\mathcal{A}$ as a $\CC$-vector spaces is given by the functions 
\begin{align*}
p_{a_1,k_1,a_2,k_2}:=\mathrm{tr}\left((\wedge^2V_1)^{\otimes a_1} \otimes \mathrm{Sym}^{k_1}(V_1) \otimes (\wedge^2V_2)^{\otimes a_2} \otimes \mathrm{Sym}^{k_2}(V_2)\right)
\end{align*}
where $a_1,a_2 \in \ZZ$, $k_1,k_2 \in \ZZ_{\geq 0}$ satisfy $2a_1+k_1=2a_2+k_2$.  It is well-known that 
\begin{align*}
\mathcal{S}^{-1}(p_{a_1,k_1,a_2,k_2})=q^{-(k_1+k_2)/2}\one_{\varpi^{a_1}\OO_B^{\times} \times \varpi^{a_2}\OO_B^\times}*\one_{\{X \in \OO_B^2: \nu(X_i)=\varpi^{k_i}\OO_F^\times\}}.
\end{align*}
Now the definition of the Satake isomorphism amounts to taking a constant term and then applying the Satake isomorphism for a split torus.  With this in mind it is not hard to see that the preimage of $\mathcal{S}^{-1}(p_{a_1,k_2,a_2,k_2})$ in $C_c^\infty(G(F)//G(\OO_F))$ is just $\mathcal{S}^{-1}(p_{a_1,k_2,a_2,k_2})|_{G(F)}$.  

We now check compatibility with the local Langlands correspondence.  Let $\pi_0$ be the unique unramified representation in $\mathrm{Res}(\widetilde{\pi})$. Let
$\varphi$ be a spherical vector in the space of $\widetilde{\pi}$.  Then $\varphi$ is in the $\pi_0$-isotypic subspace of $\pi$.  It is not hard to see that
$$
\widetilde{\pi}(\mathcal{S}^{-1}(p_{a_1,k_1,a_2,k_3}))\varphi=\pi_0(\mathcal{S}^{-1}(p_{a_1,k_1,a_2,k_3})|_{G(F)})\varphi
$$
provided that we normalize measures so that $\mathrm{meas}(\widetilde{G}(\OO_F))=\mathrm{meas}(G(\OO_F))$.  This proves compatibility of the transfer with the correspondence between unramified representations induced by $\mathrm{Res}:{}^L\widetilde{G} \to {}^LG$.  
\end{proof}

\begin{cor} \label{cor:RS:Lfuncs}Assume that $F$ is nonarchimedean and $\pi$ is an unramified irreducible admissible representation of $G(F)$.  Then $\pi$ is a transfer of an unramified representation $\widetilde{\pi}=\pi' \otimes \pi''$ of $\widetilde{G}(F)$.  One has
$$
L(s,\pi,\rho)=L(s,\pi' \times \pi'')
$$
where the $L$-function on the right is the usual Rankin-Selberg $L$-function. 
\end{cor}

\begin{proof} To prove that $\pi$ is a transfer of an unramified representation we note that every semisimple matrix in ${}^LG^{\circ}$ is the image of a semisimple matrix in 
${}^L\widetilde{G}^{\circ}$.  The last assertion is immediate from Lemma \ref{lem:compat}.
\end{proof}

Let $C_{ac}^\infty(G(F))$ denote the subspace of functions $f$ in  $C^\infty(G(F))$ such that for all compact subsets $\Omega \subset F^\times$ the restriction of $f$ to $\nu^{-1}(\Omega)$ is compactly supported.  
If $F$ is nonarchimedean, then it is easy to see that 
$$
\one_{M(\OO_F)} \in C_{ac}^\infty(G(F)).
$$

\begin{lem} \label{lem:basic:comp}
Assume $\pi$ is an unramified irreducible admissible representation of $G(F)$.  For $\mathrm{Re}(s)$ sufficiently large one has 
\begin{align*}
\mathrm{tr}\,\pi_{s+1}(\one_{M(\OO_F)})=\frac{L(s,\pi,\rho)}{L(2s,\chi_{\pi})}.
\end{align*}
\end{lem}

\begin{proof}
For $\OO_F$-algebras $R$ let
\begin{align}
\widetilde{M}(R):=\{X \in (\OO_B \otimes_{\OO_F} R)^2: \nu(X_1)=a \nu(X_2) \textrm{ for some }a \in R^\times\}.
\end{align}
The explicit description of $\mathrm{Res}^{*}$ in Lemma \ref{lem:compat} implies that 
$$
\mathrm{Res}^{*}(\one_{\widetilde{M}(\OO_F)})=\one_{M(\OO_F)}.
$$
Thus in view of Lemma \ref{lem:compat} and Corollary \ref{cor:RS:Lfuncs} it suffices to verify that for unramified representations $\pi' \otimes \pi''$ of $\widetilde{G}(F)$ one has
\begin{align} \label{for:Cauchy}
\frac{L(s,\pi' \times \pi'')}{L(2s,\chi_{\pi'}\chi_{\pi''})}=\mathrm{tr}\,(\pi' \otimes \pi'')|\omega|^{s+1}(\one_{\widetilde{M}(\OO_F)}).
\end{align}
One has
$$
\one_{\widetilde{M}(\OO_F)}=\sum_{k=0}^\infty \one_{\{X \in \OO_B^2:\nu(X_i)=\varpi^k\OO_F^\times\}}.
$$
The image of this under the Satake isomorphism is
\begin{align}
\sum_{k=0}^\infty q^k\mathrm{tr}\left(\mathrm{Sym}^k(V_1) \otimes \mathrm{Sym}^k(V_2)\right)
\end{align}
in the notation of the proof of Lemma \ref{lem:compat} so we deduce \eqref{for:Cauchy} using the Cauchy identity.
\end{proof}

Assume for the moment that $F$ is archimedean and that $B$ is split.  Let $\pi' \otimes \pi''$ 
be an irreducible unitary representation of $G(F)$.  The \textbf{analytic conductor of $L(s,\pi' \times \pi'')$} is
\begin{align} \label{AC0}
C(\pi' \times \pi'',s):=\prod_{i=1}^2\prod_{j=1}^2\left|1+\mu_{\pi_{1i} \times \pi_{2j}}+s \right|
\end{align}
where the $\mu_{\pi_{1i} \times \mu_{2j}}$ are the complex numbers such that 
$$
L(s,\pi' \times \pi'')=\prod_{i=1}^2 \prod_{j=1}^2 \zeta(s+\mu_{\pi_{1i} \times \pi_{2j}})
$$
where $\zeta$ is the Tate local $L$-function of the trivial character.  We define the analytic conductor $C(\pi,\rho,s)$ of $L(s,\pi,\rho)$ to be the analytic conductor of $L(s,\pi' \times \pi'')$ where $\pi$ is a transfer of $\pi' \otimes \pi''$:
\begin{align} \label{AC}
C(\pi,\rho,s):=C(\pi' \times \pi'',s).
\end{align}
One checks that $L(s,\pi' \times \pi'')$ depends only on $\pi$, not on $\pi' \otimes \pi''$, so this is well-defined.  For the following lemma we refer to \cite[Lemma 4.4]{GetzApproach}:

\begin{lem}\label{lem:AC:rap}
Assume that $F$ is archimedean and $B$ is split.  Then for any 
 $f \in C_c^\infty(G(F))$ that is finite under a maximal compact subgroup of $G(F)$ one has
 $$
|\mathrm{tr}\,\pi(f)| C(\pi,\rho,0)^N \ll_{f,N} 1.
 $$ \qed
\end{lem}
\noindent Technically speaking this is only proved for $\GL_2$ in loc.~cit., but the argument goes through with only minor modifications.

We also record the following trivial lemma:
\begin{lem}  \label{lem:triv} Let $v$ be an infinite place of $F$ and let $f \in C_c^\infty(G(F))$ be finite under a maximal compact subgroup of $G(F)$.  If $B_{F_v}$ is a division algebra then for all $N \geq 0$
$$
|\mathrm{tr}\,\pi(f)| C(\chi_{\pi})^N \ll_{F,N} 1
$$ 
Here $C(\chi_{\pi})$ is the analytic conductor of the central character of $\pi$.  \qed
\end{lem}

\subsection{Global $L$-packets}
\label{sec:Global}

 By restricting functions one obtains a $G(\A_F)$-equivariant map
\begin{align*}
\mathrm{res}:L^2([\widetilde{G}])^{\mathrm{sm}} \lto L^2([G])^{\mathrm{sm}}
\end{align*}
where the superscript $\mathrm{sm}$ denotes the subspace of smooth vectors.  
We say that an automorphic representation $\pi$ of $A_G \backslash G(\A_F)$ is a \textbf{transfer} of an automorphic representation $\widetilde{\pi}$ of $A_{\widetilde{G}} \backslash \widetilde{G}(\A_F)$ and write
\begin{align}
\pi \in \mathrm{res}(\widetilde{\pi})
\end{align} 
if there is a nonzero vector in the $\pi$-isotypic subspace of $L^2([G])^{\mathrm{sm}}$ that is in the image of the $\widetilde{\pi}$-isotypic subspace of $L^2([\widetilde{G}])^{\mathrm{sm}}$.  Every automorphic representation $\pi$ of $A_{G} \backslash G(\A_F)$ is is the transfer of some automorphic representation $\widetilde{\pi}$ of $A_{\widetilde{G}} \backslash \widetilde{G}(\A_F)$ \cite[Theorem 4.13]{Hiraga:Saito}.

If $\pi$ is a transfer of $\widetilde{\pi}$, then it is obvious that $\pi_v$ is a transfer of $\widetilde{\pi}_v$ for all places $v$.

\section{Spectral expansion} \label{sec:spectral}

For the rest of this paper $F$ denotes a number field (not a local field).  We also assume for the remainder of the paper that $B$, our quaternion algebra, is nonsplit over $F$.
Let $f_S \in C_c^\infty(G(F_S))$.  
We require a spectral expansion of the term
\begin{align} \label{overX}
\sum_{\gamma \in M(F)} f_S\one_{M(\widehat{\OO}_F^S)}\left(\frac{\gamma}{\Delta(\sqrt{X})} \right)=\sum_{\gamma \in G(F)}f_S\one_{M(\widehat{\OO}_F^S)}\left(\frac{\gamma}{\Delta(\sqrt{X})} \right)
\end{align}
in Theorem \ref{thm:main:asymp}. Here we are using our assumption that $f_S \in C_c^\infty(G(F_S))$ to obtain the equality.  Let $f:=f_S\one_{G(\widehat{\OO}_F^S)}$.  Note that this is different from the convention of Theorem \ref{thm:main:asymp}, in which $f$ was a Schwartz function on $B_{\A_F}^2$.  
The spectral expansion of 
$$
K_f(g_1,g_2)=\sum_{\gamma \in G(F)}f(g_1^{-1}\gamma g_2)
$$
is given by 
\begin{align} \label{cusp:exp}
K_f(g_1,g_2)=\frac{1}{2\pi i}\sum_{\pi} \int_{i\RR}K_{\pi_s(f)}(g_1,g_2)ds
\end{align}
where the sum is over isomorphism classes of automorphic representations
 $\pi$ of $A_{G} \backslash G(\A_F)$ and
 \begin{align} \label{cusp:exp2}
 K_{\pi_s(f)}(g_1,g_2)=
 \sum_{\varphi \in \mathcal{B}_\pi}\pi_s(f)\varphi(g_1)\bar{\varphi}(g_2)
 \end{align}
 where $\mathcal{B}_\pi$ is an orthonormal basis of the $\pi$-isotypic subspace of $L^2([G])$. Moreover, the sum
 \begin{align} \label{abs:conv}
\sum_{\pi} \int_{i\RR}|K_{\pi_s(f)}(g_1,g_2)|ds 
 \end{align}
 converges uniformly for $g_1,g_2$ in compact subsets of $G(\A_F)$.  This is well-known, and proved in a more general setting in \cite[Lemma 4.4]{Arthur:TFI}.
 
For $X \in \RR_{>0}$ set
$$
f_{X}(g):=f\left( \frac{g}{\Delta(\sqrt{X})}\right).
$$
Moreover for an automorphic representation $\pi$ of $G(\A_F)$ let $\chi_{\pi}$ denote its central character, so $\chi_{\pi}(z)\pi(g)=\pi(zg)$ for $(z,g) \in Z_G(\A_F) \times G(\A_F)$.

We are now in a position to state the main theorem of this section:
\begin{thm} \label{thm:basic:id} 
Assume that  $\Phi=\Phi_0 \otimes f_S\one_{\widehat{\OO}_B^S}$ satisfies the standard assumptions.  
 Let $\Omega \subset G(\A_F) \times G(\A_F)$ be a compact set.
For $\sigma>2,\varepsilon>0$, $X \in \RR_{>0}$ and $(g_1,g_2) \in \Omega$ the sum 
\begin{align}\label{spec:side}
\frac{1}{2\pi i}\sum_\pi \int_{i\RR+\sigma}\frac{L^S(s,\pi,\rho)}{L^S(2s,\chi_{\pi})}
X^{s+1}K_{\pi_{s+1}(f)}(g_1,g_2)ds
\end{align}
is equal to $O_{\varepsilon,\Phi,\Omega}(X^{3/2+\varepsilon})$ plus
\begin{align} \label{geo:side}
&X^3\left(|\omega(g_1g_2^{-1})|^4\mathcal{I}_{-1}(\Phi)(0)-|\omega(g_1g_2^{-1})|^3\mathcal{I}_{-1}(\Phi^{\mathrm{sw}})(0) \right)\\&+
X^2\left(|\omega(g_1g_2^{-1})|^4\sum_{\gamma \in M(F)}\mathcal{I}(\Phi)( g_2^{-1}\gamma g_1)-
|\omega(g_1g_2^{-1})|\sum_{\gamma \in M(F)}\mathcal{I}(\Phi^{\mathrm{sw}})(\omega(g_1^{-1}g_2)g_2^{-1}\gamma g_1)\right).
\end{align}
\end{thm}

The proof will be given after some preparation.

\begin{lem} \label{lem:sat:comp}
Let $\pi$ be an unramified irreducible unitary representation of $G(\A_F^S)$.  Then for $\mathrm{Re}(s)>2$ one has
$$
\mathrm{tr}\,\pi_{s+1}(\one_{M(\widehat{\OO}_F^S)})=\frac{L^S(s,\pi,\rho)}{L^S(2s,\chi_{\pi})}.
$$
The integral defining the (rank $1$) operator $\pi_{s+1}(\one_{M(\widehat{\OO}_F^S)})$ converges absolutely for $\mathrm{Re}(s)>2$ and the infinite sum defining the right hand side converges absolutely for $\mathrm{Re}(s)>2$.  Moreover both are bounded independently of $\pi$ in this open half plane.
\end{lem}
\noindent

\begin{proof}  
By Lemma \ref{lem:basic:comp},
 to prove the lemma it suffices to show that $\mathrm{tr}\,\pi_{s+1}(\one_{M(\widehat{\OO}_F^S)})$ converges absolutely for $\mathrm{Re}(s)>2$.  If $\pi$ is unitary and $f \in C_c^\infty(G(\A_F^S)//G(\widehat{\OO}_F^S))$ then $\mathrm{tr}\,\pi(f)$ is bounded by $\mathrm{tr}\,\pi_{\mathrm{triv}}(f)$ ($\pi_{\mathrm{triv}}$ being the trivial representation).  Thus to check the required convergence it suffices to note that
$$
\mathrm{tr}\,\pi_{\mathrm{triv},s+1}(\one_{M(\widehat{\OO}_F^S)})=\frac{\zeta^S(s+1)\zeta^S(s-1)\zeta^S(s)^2}{\zeta^S_F(2s)}
$$
for $\mathrm{Re}(s)>2$.  In fact, the sum on the right and the integral defining the operator  $\pi_{\mathrm{triv},s+1}(\one_{M(\widehat{\OO}_F^S)})$ on the left converge absolutely for $\mathrm{Re}(s)>2$.
\end{proof}

\begin{prop} \label{prop:spec:before:shift}
For $\sigma>2$ one has
$$
K_{f_S\one_{M(\widehat{\OO}_F^S)}}(g_1,g_2)=
\frac{1}{2 \pi i}\sum_\pi \int_{i\RR +\sigma }\frac{L^S(s,\pi,\rho)}{L^S(2s,\chi_{\pi})}
K_{\pi_{s+1}(f)}(g_1,g_2) ds
$$
Moreover, for $\sigma$ in the same range
$$
\sum_\pi \int_{i\RR +\sigma }\left|\frac{L^S(s,\pi,\rho)}{L^S(2s,\chi_{\pi})}
K_{\pi_{s+1}(f)}(g_1,g_2) \right|ds 
$$
converges uniformly for $(g_1,g_2)$ in compact subsets of $G(\A_F) \times G(\A_F)$.  
\end{prop}

\begin{proof}
Write $\one_{M(\widehat{\OO}_F^S)}=\sum_{m} h_m^S$ where the sum on $m$ is over $\A_F^{S \times}/\widehat{\OO}_F^{S \times}$ and 
$$
h_m^S \in C_c^\infty(G(\A_F^S)//G(\widehat{\OO}_F^S))
$$ 
is supported on the set of elements $g \in G(\A_F^S)$ such that $\omega (g) \in m$.  For each $m$ we apply \eqref{cusp:exp} and a contour shift to see that 
\begin{align*}
K_{f_Sh^S_m}(g_1,g_2)&= \frac{1}{2\pi i}\sum_\pi \int_{i\RR+\sigma}
K_{\pi_{s+1}(f_Sh_m^S)}(g_1,g_2) ds\\
&=\frac{1}{2\pi i}\sum_\pi \int_{i\RR+\sigma}\mathrm{tr}\,\pi_{s+1}(h_m^S)
K_{\pi_{s+1}(f)}(g_1,g_2)ds
\end{align*}
for any $\sigma \in \RR$, in particular for $\sigma>2$.  Thus summing over $m$ and applying Lemma 
\ref{lem:sat:comp} we deduce the identity of the proposition. The absolute convergence statement follows from Lemma \ref{lem:sat:comp} and \eqref{abs:conv}.
\end{proof}

The following lemma is an easy consequence of Theorem \ref{thm:adelic:bound}:

\begin{lem} \label{lem:unif:conv} For any $\Phi \in \mathcal{S}(\A_F^{2} \oplus B_{\A_F}^2)$ the sum 
\begin{align*}
\sum_{\gamma \in M(F)}\left|\mathcal{I}(\Phi)(g_1^{-1}\gamma g_2)\right|
\end{align*}
converges uniformly for $(g_1,g_2)$ in a compact subset of $G(\A_F) \times G(\A_F)$.  \qed
\end{lem}

With all this preparation complete we can now prove Theorem \ref{thm:basic:id}:

\begin{proof}[Proof of Theorem \ref{thm:basic:id}]
Let $h:=f_S \one_{M(\widehat{\OO}_F^S)}$.  We have
\begin{align}
\sum_{\gamma \in M(F)}h\left(\frac{g_1^{-1}\gamma g_2}{\Delta(\sqrt{X})} \right)=K_{f_{SX}\one_{M(\widehat{\OO}_F^S)}}(g_1,g_2).
\end{align}
In view of Proposition \ref{prop:spec:before:shift} for $\sigma>2$ this is equal to 
\begin{align} \label{before:shift}
K_{f_{SX}\one_{M(\widehat{\OO}_F^S)}}(g_1,g_2) &=\frac{1}{2\pi i}\sum_\pi \int_{i\RR+\sigma}\frac{L^S(s,\pi,\rho)}{L^S(2s,\chi_{\pi})}
X^{s+1}K_{\pi_{s+1}(f)}(g_1,g_2)ds
\end{align}
which is the first expression in Theorem \ref{thm:basic:id}.  As for the second expression we apply Theorem \ref{thm:main:asymp} not with
$h$ and $\Phi_0$, but with 
$$
h_{g_1,g_2}(g):=f(g_1^{-1}gg_2) \textrm{ and }\Phi_{0,g_1,g_2}(x,y):=\Phi_0( \omega(g_1^{-1}g_2)x,y).
$$  
Set $\Phi_{g_1,g_2}:=\Phi_{0,g_1,g_2} \otimes h_{g_1,g_2}$.  Since $\Phi_0\otimes h$ satisfies the standard assumptions so does $\Phi_{g_1,g_2}$.
Then for any $\varepsilon>0$,
\begin{align} \label{sum}
\sum_{\gamma \in M(F)} h\left(\frac{g_1^{-1}\gamma g_2}{\Delta(\sqrt{X})} \right)
\end{align}
is equal to $O_{\varepsilon,\Phi,g_1,g_2}(X^{\tfrac{3}{2}+\varepsilon})$ 
plus
\begin{align*}
&X^3\left(\mathcal{I}_{-1}(\Phi_{g_1,g_2})(0)-\mathcal{I}_{-1}(\Phi_{g_1,g_2}^{\mathrm{sw}})(0) \right)+X^2\sum_{\gamma \in M(F)}\left( \mathcal{I}(\Phi_{g_1,g_2})(\gamma)-\mathcal{I}(\Phi^{\mathrm{sw}}_{g_1,g_2})(\gamma)\right).
\end{align*} 
A change of variables using the definition of 
$\mathcal{I}_{-1}(\Phi)$ and $\mathcal{I}(\Phi)$ then implies the identity in the theorem. 

We are left with explaining why we can replace the implicit constant in $O_{\varepsilon,\Phi,g_1,g_2}(X^{\tfrac{3}{2}+\varepsilon})$ with one that depends only on the compact set $\Omega \subset G(\A_F) \times G(\A_F)$.  To see this we note that the basic result used in bounding the error term in Theorem \ref{thm:main:asymp} is \cite[Theorem 3.2]{GetzQuad}.  
The proof of this theorem can be easily modified to yield a version that is uniform over $(g_1,g_2)$ in the compact set $\Omega$.  More precisely, the modification necessary is just the argument used in the proof of Theorem \ref{thm:adelic:bound} to replace the bound \eqref{withgs} with \eqref{wogs}.
\end{proof}

\section{Godement-Jacquet-type zeta integrals} \label{sec:GJ}

In this section we prove Theorem \ref{thm:zeta}.  
We begin by using Theorem \ref{thm:basic:id} to deduce some of the analytic behavior of $L^S(s,\pi,\rho)$ without using Rankin-Selberg theory.
For $V \in C_c^\infty(\RR_{>0})$ let 
$$
\widetilde{V}(s):=\int_0^\infty V(x)x^{s-1}dx
$$
be its Mellin transform.
\begin{thm} \label{thm:residue} Let $V \in C_c^\infty(\RR_{>0})$.  Assume 
that $\pi$ is infinite dimensional.  One has
\begin{align*}
\frac{1}{2 \pi i }\int_{i \RR+\sigma}\frac{L^S(s,\pi,\rho)}{L^S(2s,\chi_{\pi})}X^{s+1}\widetilde{V}(s+1)ds=c_{V,\pi}X^2+O_{\varepsilon,\pi}(X^{\tfrac{3}{2}+\varepsilon})
\end{align*}
for some $c_{V,\pi} \in \RR$.
\end{thm}
This theorem is a consequence of  well-known properties of the Rankin-Selberg $L$-function. We will give a new proof based on the summation formula in Theorem \ref{thm:main:asymp}. 
The value at $c_{V,\pi}$ is nonzero for some choice of $V$ if and only if $L^S(s,\pi,\rho)$ has a pole at $s=1$.  Of course one knows via Rankin-Selberg theory that $L^S(s,\pi,\rho)$ has a pole at $s=1$ if and only if $\pi$ is in  $\mathrm{res}(\pi_0 \otimes \pi_0^\vee)$ (see \S \ref{sec:Global}) for some automorphic representation $\pi_0$ of $B_{\A_F}^\times$, but unfortunately we do not know how to give a new proof of this fact.

\begin{proof}
We assume that for $a \in \RR_{>0}$ and $g \in G(\A_F)^1$ one has
$$
f_S(\Delta(a)g)=V(a)f_S^1(g)
$$
for some $f_S^1 \in C_c^\infty(G(F_S)^1)$.  We also choose a function $\widetilde{f}_S \in C_c^\infty(G(F_S))$ such that $\int_{A_G}\widetilde{f}_S(ag)=f_S^1(g)$ for $g \in G(\A_F)^1$.  
Let $\widetilde{f}:=\widetilde{f}_S\one_{G(\widehat{\OO}_F^S)}$.

We start with the identity of Theorem \ref{thm:basic:id}.  We will multiply both sides by $\overline{K}_{\pi(\widetilde{f})}(g_1,g_2)$
and integrate along $(G(F) \backslash G(\A_F)^1)^2$.  Since the domain of integration is compact  the uniform convergence statements of 
Proposition \ref{prop:spec:before:shift} and Lemma 
\ref{lem:unif:conv} imply that we are free to bring this integral inside the other sums and integrals occurring in the identity of Theorem \ref{thm:basic:id}.

The integral of the spectral side \eqref{spec:side} times $\overline{K}_{\pi(\widetilde{f})}(g_1,g_2)$ over $(G(F) \backslash G(\A_F)^1)^2$ is
\begin{align*} 
&\int_{(G(F) \backslash G(\A_F)^1)^2}\overline{K}_{\pi(\widetilde{f})}(g_1,g_2)\frac{1}{2\pi i}\sum_{\pi'} \int_{i\RR+\sigma}\frac{L^S(s,\pi',\rho)}{L^S(2s,\chi_{\pi'})}
X^{s+1}K_{\pi_{s+1}'(f)}(g_1,g_2)ds dg_1dg_2 \nonumber\\
&=\sum_{\pi'}\frac{1}{2\pi i} \int_{i\RR+\sigma}\frac{L^S(s,\pi',\rho)}{L^S(2s,\chi_{\pi})}
X^{s+1}\widetilde{V}(s+1)ds\int_{(G(F) \backslash G(\A_F)^1)^2}K_{\pi'(\widetilde{f})}(g_1,g_2)\overline{K}_{\pi(\widetilde{f})}(g_1,g_2)dg_1dg_2\\
&=\frac{1}{2\pi i} \int_{i\RR+\sigma}\frac{L^S(s,\pi,\rho)}{L^S(2s,\chi_{\pi})}
X^{s+1}\widetilde{V}(s+1)ds \, m_{\pi} \mathrm{tr}\,\pi(\widetilde{f}^**\widetilde{f})
\end{align*}
where $\widetilde{f}^*(x):=\bar{\widetilde{f}(x^{-1})}$ and $m_{\pi}$ is the multiplicity of $\pi$ in $L^2(G(F) \backslash G(\A_F)^1)$.

On the other hand the integral of $\overline{K}_{\pi(f)}(g_1,g_2)$ times \eqref{geo:side} over $(G(F) \backslash G(\A_F)^1)^2$ is
\begin{align*}
&X^3\int_{(G(F) \backslash G(\A_F)^1)^2}\overline{K}_{\pi(f)}(g_1,g_2)
\left(|\omega(g_1g_2^{-1})|^4\mathcal{I}_{-1}(\Phi)(0)-|\omega(g_1g_2^{-1})|^3\mathcal{I}_{-1}(\Phi^{\mathrm{sw}})(0) \right)dg_1dg_2\\&+
X^2\int_{(G(F) \backslash G(\A_F)^1)^2}\overline{K}_{\pi(f)}(g_1,g_2)\\& \times \left(|\omega(g_1g_2^{-1})|^4\sum_{\gamma \in M(F)}\mathcal{I}(\Phi)(g_2^{-1}\gamma g_1)-
|\omega(g_1g_2^{-1})|\sum_{\gamma \in M(F)}\mathcal{I}(\Phi^{\mathrm{sw}})(\omega(g_1^{-1}g_2)g_2^{-1}\gamma g_1)\right) dg_1dg_2
\end{align*}
Since $\pi$ is infinite dimensional, the first term here vanishes by orthogonality.  

In Theorem \ref{thm:basic:id} take $\Omega$ to be any set containing a product of two fundamental domains for $G(F)$ acting on $G(\A_F)^1$ and take $\tfrac{1}{2}>\varepsilon>0$. We have shown that for any $\varepsilon>0$ one has
\begin{align*} 
&\frac{1}{2\pi i} \int_{i\RR+\sigma}\frac{L^S(s,\pi',\rho)}{L^S(2s,\chi_{\pi})}
X^{s+1}\widetilde{V}(s+1)ds \,m_{\pi} \mathrm{tr}\,\pi(\widetilde{f}^**\widetilde{f})+O_{\varepsilon,\Phi}(X^{\tfrac{3}{2}+\varepsilon})
\\&=X^2\int_{(G(F) \backslash G(\A_F)^1)^2}\overline{K}_{\pi(f)}(g_1,g_2)\\& \times \left(|\omega(g_1g_2^{-1})|^4\sum_{\gamma \in M(F)}\mathcal{I}(\Phi)(g_2^{-1}\gamma g_1)-
|\omega(g_1g_2^{-1})|\sum_{\gamma \in M(F)}\mathcal{I}(\Phi^{\mathrm{sw}})(\omega(g_1^{-1}g_2)g_2^{-1}\gamma g_1)\right).
\end{align*}
We can choose $\widetilde{f}$ so that $m_{\pi}\mathrm{tr}\,\pi(\widetilde{f}^**\widetilde{f}) \neq 0$.
For such an $f$ we can  take $c_{V,\pi}$ to be the coefficient of $X^2$ in this equality divided by $m_{\pi}\mathrm{tr}\,\pi(\widetilde{f}^**\widetilde{f})$.
\end{proof}

With Theorem \ref{thm:residue} in hand, to prove the functional equation of Theorem \ref{thm:zeta} we proceed as in the work of Godement and Jacquet.  Let $\pi$ be an infinite dimensional representation of $G(\A_F)$ and let $\varphi$ be a smooth form in its space.  We consider first the integral
\begin{align*}
Z_a(\mathcal{I}(\Phi),\varphi):=\int_{G(\A_F)^1} \mathcal{I}(\Phi)(a g)\varphi(g)dg
\end{align*}
for $a \in A_{G}$

For the remainder of this section we choose $f_S \in C_c^\infty(G(F_S))$ such that 
for $a \in \RR_{>0}$ and $g \in G(\A_F)^1$ one has
\begin{align} \label{h:prod}
f_S(\Delta(a)g)=V(a)f_S^1(g)
\end{align}
for some $f_S^1 \in C_c^\infty(G(F_S)^1)$ and some $V \in C_c^\infty(\RR_{>0})$.  We also choose a function $\widetilde{f}_S \in C_c^\infty(G(F_S))$ such that $\int_{A_G}\widetilde{f}_S(ag)=f_S^1(g)$ for $g \in G(\A_F)^1$.  
Let $f:=f_S\one_{G(\widehat{\OO}_F^S)}$, $\widetilde{f}:=\widetilde{f}_S\one_{G(\widehat{\OO}_F^S)}$ and
 let
 $$
 \Phi:=\Phi_0 \otimes f_S\one_{\widehat{\OO}_B^S} \in \mathcal{S}(\A_F^2 \times B_{\A_F}^2)
 $$
 satisfy the standard assumptions.
 
 Write
\begin{align}
c_{V,\pi}:=\lim_{X \to \infty} \frac{1}{2\pi i X^2}\int_{i \RR+\sigma}\frac{L^S(s,\pi,\rho)}{L^S(2s,\chi_{\pi})}X^{s+1}\widetilde{V}(s+1)ds.
\end{align}
This quantity is well-defined by Theorem \ref{thm:residue}.  Of course, if we apply Rankin-Selberg theory we see that it is 
$$
\frac{\widetilde{V}(1)\mathrm{Res}_{s=1}L^S(s,\pi,\rho)}{L^S(2,\chi_{\pi})},
$$
but part of the point of the current section is to see how much we can derive directly from Theorem \ref{thm:basic:id} without assuming this.

\begin{prop} \label{prop:first:func}
One has 
\begin{align*}
Z_a(\mathcal{I}(\Phi),\varphi)=|\omega(a)|^{-3}Z_{a^{-1}}(\mathcal{I}(\Phi^{\mathrm{sw}}),\varphi^{\vee})+|\omega(a)|^{-2}c_{V,\pi}\int_{G(\A_F)^1}\widetilde{f}(g^{-1})\varphi(g)dg.
\end{align*}
\end{prop}

\begin{proof}
We take the identity of Theorem \ref{thm:basic:id} for $(g_1,g_2)=(ag,I)$ and integrate it times $\varphi(g)$ along $G(F) \backslash G(\A_F)^1$.  Since $\pi$ is an infinite-dimensional representation the coefficients of $X^3$ vanish identically and we arrive at the equality of
\begin{align} \nonumber
&\lim_{X \to \infty}\frac{1}{X^2}\int_{G(F) \backslash G(\A_F)^1}\frac{1}{2\pi i }\int_{i \RR+\sigma}\sum_{\pi'}\frac{L^S(s,\pi,\rho)}{L^S(2s,\chi_{\pi})}X^{s+1}K_{\pi_{s+1}(f)}(ag,I)\varphi(g)ds dg\\
&=\lim_{X \to \infty}\frac{1}{2\pi i X^2 }\int_{i \RR+\sigma}\frac{L^S(s,\pi^{\vee},\rho)}{L^S(2s,\chi_{\pi^\vee})}X^{s+1}|\omega(a)|^{s+1}\widetilde{V}(s+1)ds \int_{G(\A_F)^1}\widetilde{f}(g^{-1})\varphi(g)dg \nonumber \\
&=|\omega(a)|^{2}c_{V,\pi^\vee} \int_{G(\A_F)^1}\widetilde{f}(g^{-1})\varphi(g)dg
 \label{spec2}
\end{align}
and
\begin{align} \label{geo2}
\begin{split}|\omega (a)|^4
&\int_{G(F) \backslash G(\A_F)^1}\sum_{\gamma \in M(F)}\mathcal{I}(\Phi)(\gamma ag)\varphi(g)dg\\
&-|\omega(a)|\int_{G(F) \backslash G(\A_F)^1}\sum_{\gamma \in M(F)}\mathcal{I}(\Phi^{\mathrm{sw}})(\omega(ag)^{-1}\gamma ag)\varphi(g)dg.
\end{split}
\end{align}
Here in \eqref{spec2} we have moved the integral over $G(F) \backslash G(\A_F)^1$ inside the integral over $i\RR+\sigma$ and the sum over $\pi$; this is justified by the uniform convergence statement of  Proposition \ref{prop:spec:before:shift}.  

Since $B$ is a division algebra 
$$
M(F)=G(F) \amalg (0,0)
$$
is a decomposition of $M(F)$ into its $G(F)$ orbits.  The contribution of the term $(0,0)$ to both the sums in \eqref{geo2} vanishes since $\pi$ is infinite dimensional.  Unfolding the remaining terms is justified by the uniform convergence statement of Lemma \ref{lem:unif:conv} and we see that \eqref{geo2} is equal to
\begin{align}
\begin{split}
|\omega(a)|^4&\int_{ G(\A_F)^1}\mathcal{I}(\Phi)(ag)\varphi(g)dg\\
&-|\omega(a)|\int_{G(\A_F)^1}\mathcal{I}(\Phi^{\mathrm{sw}})(\omega(ag)^{-1}ag)\varphi(g)dg.
\end{split}
\end{align}
The first summand here is $|\omega(a)|^4Z_a(\mathcal{I}(\Phi),\varphi)$.  Taking a change of variable $g \mapsto g^{\vee}$ we see that the second integral is 
$|\omega(a)| Z_{a^{-1}}(\mathcal{I}(\Phi),\varphi^{\vee})$.  
\end{proof}

\begin{lem} \label{lem:abs:conv01} 
For all $\Phi \in \mathcal{S}(\A_F^2 \oplus B_{\A_F}^2)$ (not necessarily satisfying the standard assumptions) and all $s \in \CC$ the integral
\begin{align*} 
\int_{|\omega(a)| > 1}\int_{G(\A_F)^1}|\mathcal{I}(\Phi)(ag)||\omega(a)|^sda dg
\end{align*}
is bounded.
\end{lem}

\begin{proof}
Arguing as in the proof of 
Theorem \ref{thm:adelic:bound} we see that there is an $A>0$ such that for any compact set $\Omega \subset G(\A_F)$ one has a bound
\begin{align} \label{withg}
|\mathcal{I}(\Phi)(a\gamma g)| \ll_{\Phi,N,\Omega} 
|\omega(a)|^{-N/2}\prod_{v|\infty}(| \gamma |_v,1)^{-N}
\end{align}
 for any $N \in \ZZ_{ \geq 0}$ provided that $g \in \Omega$ and $a \in A_G$ satisfies $|\omega(a)| \geq 1$.  Moreover, $\mathcal{I}(\Phi)(a\gamma g)$ vanishes unless $\gamma \in \beta^{-1}\OO_B$ for some $\beta \in F^\times$ that depends only on $\Omega$.

Taking $\Omega$ to be a compact measurable fundamental domain for the action of $G(F)$ on $G(\A_F)^1$ we have
\begin{align*}
&\int_{|\omega(a)| > 1}\int_{G(\A_F)^1}|\mathcal{I}(\Phi)(ag)||\omega(a)|^sda dg\\
&=\int_{|\omega(a)| > 1}\int_{G(F) \backslash G(\A_F)^1}\sum_{\gamma \in G(F)}|\mathcal{I}(\Phi)(a\gamma g)||\omega(a)|^sda dg\\
& \ll_{N} \int_{|\omega(a)| > 1} |\omega(a)|^{s-N/2}da \sum_{\gamma \in \beta^{-1}\OO_B}\prod_{v|\infty}\max(1,|\gamma|_v)^{-N}.
\end{align*}
For any $s$ we can choose $N$ large enough that this converges.
\end{proof}

We now prove Theorem \ref{thm:zeta}:

\begin{proof}[Proof of Theorem \ref{thm:zeta}]
One has
\begin{align*} 
Z(\mathcal{I}(\Phi),s,\varphi)&=\int_{A_G}Z_a(\mathcal{I}(\Phi),\varphi)|\omega(a)|^{s+1}da\\&=\int_{|\omega(a)|< 1}Z_a(\mathcal{I}(\Phi),\varphi)|\omega(a)|^{s+1}da+\int_{|\omega(a)| > 1} Z_a(\mathcal{I}(\Phi),\varphi)|\omega(a)|^{s+1}da.
\end{align*}
The latter integral converges absolutely for all $s$ by Lemma \ref{lem:abs:conv01}.  As for the former integral by Proposition \ref{prop:first:func} we have 
\begin{align*}
&\int_{|\omega(a)|< 1}Z_a(\mathcal{I}(\Phi),\varphi)|\omega(a)|^{s+1}da\\&=\int_{|\omega(a)|< 1}Z_{a^{-1}}(\mathcal{I}(\Phi^\mathrm{sw}),\varphi^{\vee})|\omega(a)|^{s-2}da+\int_{|\omega(a)|<1}|\omega(a)|^{s-1}c_{f,\varphi}da
\end{align*}
where 
$$
c_{f,\varphi}:=c_{V,\pi}\int_{G(\A_F)^1}\widetilde{f}(g^{-1})\varphi(g)dg.
$$
For $s>1$ the last term here is $\frac{c_{f,\varphi}}{2(s-1)}$, and hence by analytic continuation it takes this value for all $s \neq 1$.  The former term is
\begin{align*}
&\int_{|\omega(a)|>1}Z_{a}(\mathcal{I}(\Phi^\mathrm{sw}),\varphi^{\vee})|\omega(a)|^{2-s}da.
\end{align*}
Since $\varphi$ is bounded  this converges absolutely for all $s$ by Lemma \ref{lem:abs:conv01}.  Thus we have obtained the meromorphic continuation of $Z(\mathcal{I}(\Phi),s,\varphi)$ with poles as specified in the statement of the theorem.  Combining the equalities above we have also obtained
\begin{align} \label{1}\begin{split}
&Z(\mathcal{I}(\Phi),s,\varphi)\\
&=\int_{|\omega(a)|>1}Z_{a}(\mathcal{I}(\Phi^\mathrm{sw}),\varphi^{\vee})|\omega(a)|^{2-s}da+\int_{|\omega(a)| > 1} Z_a(\mathcal{I}(\Phi),\varphi)|\omega(a)|^{s+1}da+\frac{c_{f,\varphi}}{2(s-1)}.\end{split}
\end{align}

In the traditional argument for the functional equation of zeta functions going back to Tate's thesis one would now argue that one is done by symmetry.  However, this will not work for us because, as already observed, if $\Phi$ satisfies the standard assumptions then $\Phi^{\mathrm{sw}}$ does not.  Moreover, the expression for the residues we obtained above is asymmetric in $\Phi$ and $\Phi^{\mathrm{sw}}$.

Instead we start over, this time with $\Phi$ replaced by $\Phi^{\mathrm{sw}}$.  We have
\begin{align*}
Z(\Phi^{\mathrm{sw}},1-s,\varphi^{\vee})&=\int_{A_G}Z_a(\mathcal{I}(\Phi^{\mathrm{sw}}),\varphi^{\vee})|\omega(a)|^{2-s}da\\&=\int_{|\omega(a)|< 1}Z_a(\mathcal{I}(\Phi^{\mathrm{sw}}),\varphi^\vee)|\omega(a)|^{2-s}da+\int_{|\omega(a)| > 1} Z_a(\mathcal{I}(\Phi^{\mathrm{sw}}),\varphi^\vee)|\omega(a)|^{2-s}da.
\end{align*}
The second summand here converges absolutely for all $s$ by Lemma \ref{lem:abs:conv01}.  By Proposition \ref{prop:first:func} the first summand is equal to 
\begin{align*}
&\int_{|\omega(a)|< 1}Z_{a^{-1}}(\mathcal{I}(\Phi),\varphi)|\omega(a)|^{-1-s}da-\int_{\omega(a)<1}|\omega(a)|^{1-s}c_{f,\varphi}da\\
&=\int_{|\omega(a)|>1}Z_{a}(\mathcal{I}(\Phi),\varphi)|\omega(a)|^{1+s}da+\frac{c_{f,\varphi}}{2(s-1)}.
\end{align*}
Here the second summand is only defined for $\mathrm{Re}(s)<1$ at first, but after evaluating the integral one obtains a function that obviously extends meromorphically to the plane.  The first summand here converges absolutely for all $s$ by Lemma \ref{lem:abs:conv01}, so we have obtained the analytic continuation of $Z(\mathcal{I}(\Phi^{\mathrm{sw}}),\varphi^{\vee},s)$ to a meromorphic function of $s$ that is holomorphic except for a simple
pole at $s=0$ with residue $-c_{\Phi,\pi}$.
Combining the previous equalities also yields
\begin{align*}
&Z(\mathcal{I}(\Phi^{\mathrm{sw}}),1-s,\varphi^{\vee})\\&=\int_{|\omega(a)|>1}Z_{a}(\mathcal{I}(\Phi),\varphi)|\omega(a)|^{1+s}da+\int_{|\omega(a)| > 1} Z_a(\mathcal{I}(\Phi^{\mathrm{sw}}),\varphi^\vee)|\omega(a)|^{2-s}da+\frac{c_{f,\varphi}}{2(s-1)}
\end{align*}
which is equal to \eqref{1}, proving the desired functional equation.
\end{proof}

\section{Four variable kernel functions}
\label{sec:four:var}
From this section onward we allow ourselves to use the entirety of Rankin-Selberg theory.  Using it we can give the following more precise 
version of Theorem \ref{thm:basic:id}. 
 Let $f_S \in C_c^\infty(G(F_S))$ and let
$$
\Phi=\Phi_0 \otimes f_S \one_{\widehat{\OO}_F^S}
$$
satisfy the standard assumptions.  Recall that $S$ is assumed to satisfy \textbf{A($S$)} from \S \ref{ssec:zeta}

\begin{thm} \label{thm:body:PS} Assume that $f_S$ is finite under a maximal compact subgroup of $G(F_S)$.  One has
\begin{align}
\label{spec:side2}&\sum_{\pi} \frac{\mathrm{Res}_{s=1}L^S(s,\pi,\rho)}{\zeta^S(2)}
K_{\pi_{2}(f)}(g_1,g_2)\\\label{geo:side2}&=
|\omega(g_1g_2^{-1})|^4\sum_{\gamma \in M(F)}\mathcal{I}(\Phi)(g_1\gamma g_2^{-1})-
|\omega(g_1g_2^{-1})|\sum_{\gamma \in M(F)}\mathcal{I}(\Phi^{\mathrm{sw}})(\omega(g_1^{-1}g_2)g_1\gamma g_2^{-1})
\end{align}
where the sum on $\pi$ is over isomorphism classes of automorphic representations of $A_{G} \backslash G(\A_F)$
that are irreducible constituents of
$\pi_0 \otimes \pi_0^{\vee}|_{G(\A_F)}$ for some automorphic representation $\pi_0$ of $ B_{\A_F}^\times$.
\end{thm}

Note that in the theorem $K_{\pi_2(f)}(g_1,g_2)$ is $K_{\pi_s(f)}(g_1,g_2)$ evaluated at $s=2$.  This notation also occurs below.
We view this as a geometric expansion of a four-variable kernel function.  
To explain this, note that each kernel $K_{\pi_2(f)}(g_1,g_2)$
is controlled by a single automorphic representation of $B_{\A_F}^\times$, namely $\pi_0$, but, up to center, there are four copies of $B_{\A_F}^\times$ that can be integrated over in the kernel. These extra copies will be put to good use in \S \ref{sec:nabel} below.

The main step in deriving Theorem \ref{thm:body:PS} from Theorem \ref{thm:basic:id} is the following proposition:

\begin{prop} \label{prop:ds}  Assume that $f_S$ is finite under a maximal compact subgroup of $G(F_S)$.  For any $\varepsilon>0$ the sum \eqref{spec:side} is equal to $O_{\varepsilon}(X^{\tfrac{3}{2}+\varepsilon})$ plus
\begin{align*}
&\sum_{\pi} \left(\frac{\mathrm{Res}_{s=2}L^S(s,\pi,\rho)}{\zeta^S(4)}
X^{3}K_{\pi_3(f)}(g_1,g_2)+\frac{\mathrm{Res}_{s=1}L^S(s,\pi,\rho)}{\zeta^S(2)}
X^{2}K_{\pi_2(f)}(g_1,g_2)\right)
\end{align*}
where the sum on $\pi$ is as in Theorem \ref{thm:body:PS}.  The sum is absolutely uniformly convergent for $(g_1,g_2)$ in compact subsets of $G(\A_F) \times G(\A_F)$.
\end{prop}

Indeed, assuming this proposition, Theorem \ref{thm:body:PS} follows upon comparing the coefficients of $X^2$ in Theorem \ref{thm:basic:id} and Proposition \ref{prop:ds}.

\begin{proof}
Assume that $\pi$ is a subrepresentation of $\pi' \otimes \pi''|_{G(\A_F)}$.  If $\pi$ is infinite dimensional (which is to say that at least one of $\pi'$ and $\pi''$ is infinite dimensional) then the Rankin-Selberg $L$-function $L^S(s,\pi,\rho)$ is holomorphic in the plane except for a possible  simple pole at $s=1$ which occurs if and only if $\pi' \otimes \pi'' \cong \pi_0  \otimes \pi_0^{\vee}$.  If $\pi$ is finite dimensional, then $\pi=\chi \circ \omega$ for some character $\chi \in \widehat{[\GG_m]}$.  In this case
\begin{align} \label{Rankin:split}
L^S(s,\pi,\rho)=L^S(s+1,\chi)L^S(s,\chi)^2L^S(s-1,\chi)
\end{align}
which again is holomorphic in the plane except for possible poles at $s=0,1,2$ which can only occur if $\chi =1$.  Since $S$ contains finite places, $L^S(s-1,\chi)$ vanishes at $s=1$, so the pole of $L^S(s,\pi,\rho)$ at $s=1$ is at worst simple.

These comments on the residues of $L^S(s,\pi,\rho)K_{\pi_{s+1}(f)}(g_1,g_2)$ together with a contour shift imply that 
\eqref{spec:side} is equal to the sum in the statement of the proposition plus
\begin{align} \label{after:shift}
\frac{1}{2\pi i}\sum_\pi \int_{i\RR+\tfrac{1}{2}+\varepsilon}\frac{L^S(s,\pi,\rho)}{L^S(2s,\chi_{\pi})}
X^{s+1}K_{\pi_{s+1}(f)}(g_1,g_2) ds .
\end{align}
We will show that this is $O_{\varepsilon}(X^{\tfrac{3}{2}+\varepsilon})$.  We will first give the argument when $B_{F_v}$ splits for all $v|\infty$, and then explain how to alter it when we remove this assumption.

For $\pi$ an automorphic representation of $A_G \backslash G(\A_F))$ let
\begin{align*}
C(\pi_\infty,\rho):=\prod_{v|\infty}C(\pi_v,\rho) 
\end{align*}
where $C(\pi_v,\rho)$ is the analytic conductor of $L(s,\pi_v,\rho)$ as in \eqref{AC}.

For $\pi$ occurring in the restriction of $\pi' \otimes \pi''$ (an automorphic representation of $(B_{\A_F}^\times)^2$)
we claim that one has a preconvexity bound
$$
(s-1)^{\delta_{1,\pi',\pi''^{\vee}}}(s-2)^{\delta_{2,\pi',\pi''^{\vee}}}L^S(s,\pi,\rho) \ll_f C(\pi_\infty,\rho)^\delta
$$
valid for sufficiently large $\delta>0$ and $\tfrac{1}{2} \leq \mathrm{Re}(s) \leq 2$, where 
\begin{align*} 
\delta_{1,\pi',\pi''^\vee}&=\begin{cases}1 & \textrm{ if }\pi' \cong \pi''^\vee \textrm{ and }\pi' \textrm{ is infinite dimensional}\\ 0 &\textrm{ otherwise.} \end{cases}\\
\delta_{2,\pi',\pi''^\vee}:&=\begin{cases} 1 & \textrm{ if }\pi' \cong \pi''^{\vee} \textrm{ and }\pi' \textrm{ is finite dimensional}\\
0 &\textrm{ otherwise.} \end{cases}
\end{align*}
We remind the reader that $S$ contains finite places, so the pole at $s=1$ is only simple, and not of order $2$, when $\pi'$ and $\pi''$ are finite dimensional.  We also note that we are only claiming the bound for $\pi$ contributing to our sum; this set of $\pi$ has the property that the nonarchimedean parts of their conductors are bounded in terms of our test function $f$ (this is why we have put in a subscript $f$ in the bound). 

 If $\pi'$ and $\pi''$ are both infinite dimensional then this is trivial for $\mathrm{Re}(s)>1$ and for $\tfrac{1}{2} \leq \mathrm{Re}(s) \leq 1$ this is proven in \cite[\S 1]{BrumleyNarrow}. If exactly one of $\pi'$ and $\pi''$ is infinite dimensional, say $\pi'$, and $\pi''=\chi \circ \nu$, then 
$$
L^S(s,\pi,\rho)=L^S(s+\tfrac{1}{2},\pi' \otimes \chi)L^S(s-\tfrac{1}{2},\pi' \otimes \chi)
$$
and we can again use the same reference.  If $\pi=\chi \circ \omega$  then we have the identity \eqref{Rankin:split}. We can then use 
 \cite[\S III.6, Theorem 14A]{Moreno:Advanced:analytic:number:theory}.

Thus by dominated convergence, to complete the proof that \eqref{after:shift} is $O_{\varepsilon}(X^{\tfrac{3}{2}+\varepsilon})$ it suffices to show that for any 
$f \in C_c^\infty(A_{G} \backslash G(\A_F))$ one has that
\begin{align*}
\sum_{\pi}C(\pi_\infty,\rho)^N |K_{\pi(f)}(g_1,g_2)|
\end{align*}
is bounded for any $N>0$.  This will also imply the absolute uniform convergence of the sum over residues if we can obtain a bound uniform for $(g_1,g_2)$ in a compact set.

By a standard argument (compare the proof of \cite[Theorem 3.1]{GetzHahnSRTF}), to prove this 
it suffices to show that for any $N>0$ and $h \in C_c^\infty(A_{G} \backslash G(\A_F))$ the sum
\begin{align} \label{is:bounded}
\sum_{\pi=\pi' \otimes \pi''} C(\pi,\rho)^N K_{\pi(f*f^*)}(g,g) 
\end{align} 
is bounded, uniformly for $g$ in a compact set.
Here
$h^*(g):=\bar{h}(g^{-1})$, the bar denoting complex conjugation.  The Casimir eigenvalue of $\pi$ is bounded by a constant times a power of $C(\pi,\rho)$.  This was proven for $\GL_n$ in \cite[Lemma 4.5]{GetzApproach}, and can be proven in the current setting by a trivial modification of the argument.  In view of the Weyl law we see that applying \cite[(15')]{Godement_cusp} (stated in adelic language in \cite[Theorem 3.5]{GetzHahnSRTF}), to prove the boundedness of \eqref{is:bounded} it suffices to show that $\mathrm{tr}\,\pi(f*f^*)$ is rapidly decreasing as a function of $C(\pi_\infty,\rho)$. This is the content of Lemma \ref{lem:AC:rap}.  

Now to complete the proof we explain how to modify the argument when $B$ is nonsplit at an infinite place.  Let $S \subset \infty$ be the maximal subset such that $B_{F_v}$ is nonsplit for all $v \in S$.  Then preconvexity bound is given in terms of the Jacquet-Langlands transfer of an extension of $\pi_\infty$ to $(B_{F_\infty}^{\times})^2$. On the other hand, our analytic control on  $\mathrm{tr}\,\pi_\infty(f*f^*)$ comes from Lemma \ref{lem:triv}, which is a statement about $\pi_\infty$ itself, and not its transfer. 

To overcome this difficulty we recall we have assumed that $f$ is finite under a maximal compact subgroup of $G(F_S)$.  Thus there  
is a finite set $\mathcal{R}$ of irreducible unitary representations of $G(F_S)$ such that any $\pi$ contributing to our sum has the property that a unitary twist of $\pi_S$ is isomorphic to a representation in $\mathcal{R}$.  Twisting by characters is compatible with the Jacquet Langlands correspondence in the natural way.  Moreover one has that
$$
\mathrm{tr}\pi_{Sit}(f) 
$$
is rapidly decreasing as a function of $t \in \RR$ for any fixed $\pi_S$ and $f \in C_c^\infty(G(F_S))$, 
and 
$$
C((\pi'\times \pi'',s)\otimes |\nu|^{it})=C((\pi'\times \pi'',2it+s).
$$
Using these observations it is not hard to modify the argument above.
\end{proof}

\section{A nonabelian trace formula} \label{sec:nabel}
    
In this section we prove Theorem \ref{thm:ntf}.  We place ourselves in the setting of \S \ref{ssec:ntf}.  Thus we assume that there is a subfield $k \leq F$ such that $F/k$ is Galois with Galois group 
$$
\Gal(F/k)=\langle\iota, \tau \rangle;
$$ 
that is, $\Gal(F/k)$ is generated by two elements.    Assume moreover that $B_1$ is a division algebra over $k$ such that $B:=B_1 \otimes_kF$ is nonsplit (i.e.~again a division algebra).  We define $G_0$, $\theta:G_0 \to G_0$, and the action of $G_0$ on $G$ via $\theta$-conjugation as before.
For $\gamma \in G(F)=\mathrm{Res}_{F/k}G(k)$ we let $G_{0\gamma}$ be the stabilizer of $\gamma$ under this action.  Since $B$ is a division algebra, a standard argument implies that $G_{0\gamma}$ is reductive and anisotropic modulo center.

For suitable smooth test functions $f$ on $G(\A_F)=\mathrm{Res}_{F/k}G(\A_k)$ we let
\begin{align} \label{TO}
\mathrm{TO}_{\gamma}(f):=\int_{G_{0\gamma}(\A_k) \backslash G_0(\A_k)} f(g^{-1}\gamma g^{\theta})d\dot{g}
\end{align}
be the usual twisted orbital integral.  In addition to depending on the choice of a Haar measure on $G_0(\A_k)$
it depends on a choice of Haar measure $dt_\gamma$ on $G_{0\gamma}(\A_k)$. 

\begin{lem} \label{lem:abs:conv}
For $\Phi \in \mathcal{S}(\A_F^2 \oplus B_{\A_F}^2)$ one has
\begin{align*}
\sum_{\gamma} \mathrm{meas}([G_\gamma])\int_{G_{0\gamma}(\A_k) \backslash G_0(\A_k)}|\mathcal{I}(\Phi)(g^{-1}\gamma \theta(g))|d\dot{g}<\infty
\end{align*}
where the sum on $\gamma$ is over a set of representatives for the $G_0(F)$ orbits in $G(F)$ under $\theta$-conjugation.
In particular $\mathrm{TO}_{\gamma}(\mathcal{I}(\Phi))$ is well-defined.
\end{lem}

\begin{proof}
By unfolding we see that
\begin{align*}
\sum_{\gamma} \mathrm{meas}([G_\gamma])\int_{G_{0\gamma}(\A_k) \backslash G_0(\A_k)}|\mathcal{I}(\Phi)(g^{-1}\gamma \theta(g))|d\dot{g}=
\int_{[G_0]}\sum_{\gamma \in G_0(k)}|\mathcal{I}(\Phi)(g^{-1}\gamma \theta(g))|dg.
\end{align*}
The latter integral is absolutely convergent by Lemma \ref{lem:unif:conv} because $[G_0]$ is compact.
\end{proof}

Recall the definition of the nonabelian trace \eqref{ntr} from \S \ref{ssec:ntf}
Let 
$$
\mathrm{S}G(R):=\{ g \in B_R^{\times}: \nu(g)=1\}.
$$
 To analyze the nonabelian trace it is convenient to first state a lemma on restrictions of representations of $B_{\A_F}^\times$ to $\mathrm{S}G(\A_F)$:
 \begin{lem} \label{lem:restri}
 Let $\pi$ be an automorphic representation of $A_{\GG_m} \backslash B_{\A_F}^\times$.  Its restriction to $\mathrm{S}G(\A_F)$ is a direct sum of admissible representations of $\mathrm{S}G(\A_F)$.  Two representations 
 $\pi'$ and $\pi''$ of $B_{\A_F}^\times$ have a common constituent when restricted to $\mathrm{S}G(\A_F)$ if and only if $\pi' \cong \pi'' \otimes \chi$ for some character $\chi \in \widehat{[\GG_{mF}]}$.  
 \end{lem}

\begin{proof}
For the assertion that the restriction of the representation to $\mathrm{S}G(\A_F)$ breaks into a direct sum
 see \S \cite[Chapter 2]{Hiraga:Saito}.   Moreover, in the same reference it is proven that if the restrictions of $\pi'$ and $\pi''$ to $\mathrm{S}G(\A_F)$ have a common constituent, then $\pi' \cong \pi'' \otimes \chi$ for some character $\chi:\A_F^\times \to \CC^\times$ (they do not prove that this character is invariant under $F^\times$).    To prove that $\chi$ must be invariant under $A_{\GG_m}F^\times$ we apply \cite[Theorem 4.1.2]{Ramakrishnan:RS} and the Jacquet-Langlands correspondence.
\end{proof}

For $\pi$ an admissible representation of $B_{\A_F}$ and $\xi \in \Gal(F/k)$ let 
$$
\pi^{\xi}(g):=\pi(\xi(g))
$$
be its Galois conjugate.  
The key property of the nonabelian trace is that it is nonzero only for $\pi \in \mathrm{res}(\pi_0 \otimes \pi_0^{\vee})$ where $\pi_0$ is isomorphic to its conjugates under $\Gal(F/k)$ up to a twist by a Hecke character:
\begin{lem} \label{lem:ntr} Let $\pi \in \mathrm{res}(\pi_0 \otimes \pi_0^{\vee})$ where $\pi_0$ is an automorphic representation of $B_{\A_F}^\times$.
If 
\begin{align}
\mathrm{ntr}\,\pi(f) \neq 0
\end{align}
for some $f$ then 
\begin{align*}
\pi_0 \cong \pi_0^{\iota} \otimes \chi_1 \cong \pi_0^{\tau} \otimes \chi_2
\end{align*}
for some $\chi_1,\chi_2 \in \widehat{[\GG_{mF}]}$ such that $(\chi_1\chi_2^{-1})^2|_{[\GG_{mk}]}=1$.
\end{lem}
Here we have used the notation of \S \ref{sec:Global}.

\begin{proof}
Let $V_{\pi} \subset L^2(G(F) \backslash G(\A_F)^1)$ be the space of $\pi$ and let 
$$
\bar{V}_{\pi}^{\theta}=\{  \overline{\varphi}^{\theta} \in L^2(G(F) \backslash G(\A_F)^1):\varphi \in V_{\pi} \}
$$  
where $\overline{\varphi}^{\theta}(x):=\overline{\varphi}(\theta(x))$.
The space $\overline{V}^{\tau}_{\pi}$ is a model of the representation $\pi^{\vee \theta}$ if we 
let $(G(\A_F)^2)^1$ act on the space via the regular action:
$$
R(g)\overline{\varphi}^{\theta}(x):= \overline{\varphi}^{\theta}(xg).
$$
One has a $\CC$-bilinear map
\begin{align*}
V_{\pi}^{\mathrm{sm}} \times (\overline{V}_{\pi}^{\theta})^{\mathrm{sm}} &\lto \CC\\
( \varphi_1,\overline{\varphi}_2^{\theta}) & \longmapsto \int_{[G_0]}\varphi_1(g)\overline{\varphi}_2(\theta(g))dg.
\end{align*}
 This pairing is $G_0(\A_k)$-invariant.
In particular, if $\mathrm{ntr}\,\pi(f)$ is nonzero then there is a $G_0(\A_k)$-invariant linear form on the space of $\pi \otimes \pi^{\vee \theta}$.  

Now $\pi|_{G_0(\A_k)}$ and $\pi^{\theta}|_{G_0(\A_k)}$ decompose into finite sums of irreducible representations of $G_0(\A_k)$ since this is even true of the restrictions to $\mathrm{S}G(\A_F) \times \mathrm{S}G(\A_F)$ by Lemma \ref{lem:restri}.  If $\mathrm{ntr}\,\pi(f) \neq 0$ then 
$
\pi_0 \otimes \pi_0^{\vee}|_{G_0(\A_k)}$ 
and $\pi_0^{\iota} \otimes \pi_0^{\vee \tau}|_{G_0(\A_k)}$ have a constituent in common.  Thus the lemma follows from Lemma \ref{lem:restri}.
\end{proof}

We now prove Theorem \ref{thm:ntf}:

\begin{proof}[Proof of Theorem \ref{thm:ntf}]
In the proof of Lemma \ref{lem:abs:conv} we proved that 
\begin{align} \label{geo:conv3}
\int_{[G_0]}\sum_{\gamma \in G_0(F)}|\mathcal{I}(\Phi)(g^{-1}\gamma \theta(g))|dg
\end{align}
is finite.
In addition, since $[G_0]$ is compact we can apply Proposition \ref{prop:ds} to see that 
\begin{align} \label{spec:conv3}
\sum_{\pi} \left|\frac{\mathrm{Res}_{s=1}L^S(s,\pi,\rho)}{\zeta^S(2)}\right|
\int_{[G_0]}\left|K_{\pi_2(f)}(g,\theta(g))\right|dg
\end{align}
is finite.

Now take the identity of Theorem \ref{thm:body:PS} and integrate it over
$$
\{(g,\theta(g)):g \in G_0(\A_k)\}.
$$
Since \eqref{geo:conv3} and \eqref{spec:conv3} are finite we can bring the integral over $[G_0]$ inside the sums in \eqref{geo:side2} and \eqref{spec:side2} and deduce the theorem.  
\end{proof}


\bibliography{../refs}{}
\bibliographystyle{alpha}

\end{document}